\pdfoutput=1
\RequirePackage{ifpdf}
\ifpdf 
\documentclass[pdftex]{sigma}
\else
\documentclass{sigma}
\fi

\numberwithin{equation}{section}

\newtheorem{Theorem}{Theorem}[section]
\newtheorem{Corollary}[Theorem]{Corollary}
 { \theoremstyle{definition}
\newtheorem{Remark}[Theorem]{Remark} }

\begin{document}

\allowdisplaybreaks

\newcommand{\arXivNumber}{1606.08352}

\renewcommand{\thefootnote}{}

\renewcommand{\PaperNumber}{030}

\FirstPageHeading

\ShortArticleName{GKZ Hypergeometric Series for the Hesse Pencil}

\ArticleName{GKZ Hypergeometric Series for the Hesse Pencil,\\ Chain Integrals and Orbifold Singularities\footnote{This paper is a~contribution to the Special Issue on Modular Forms and String Theory in honor of Noriko Yui. The full collection is available at \href{http://www.emis.de/journals/SIGMA/modular-forms.html}{http://www.emis.de/journals/SIGMA/modular-forms.html}}}

\Author{Jie ZHOU}

\AuthorNameForHeading{J.~Zhou}

\Address{Perimeter Institute for Theoretical Physics, Waterloo, Ontario N2L 2Y5, Canada}
\Email{\href{mailto:jzhou@perimeterinstitute.ca}{jzhou@perimeterinstitute.ca}}

\ArticleDates{Received October 01, 2016, in f\/inal form May 14, 2017; Published online May 20, 2017}

\Abstract{The GKZ system for the Hesse pencil of elliptic curves has more solutions than the period integrals. In this work we give dif\/ferent realizations and interpretations of the extra solution, in terms of oscillating integral, Eichler integral, chain integral on the elliptic curve, limit of a period of a~certain compact Calabi--Yau threefold geometry, etc. We also highlight the role played by the orbifold singularity on the moduli space and its relation to the GKZ system.}


\Classification{14J33; 14Q05; 30F30; 34M35}

\renewcommand{\thefootnote}{\arabic{footnote}}
\setcounter{footnote}{0}

\section{Introduction}\label{secintro}

The GKZ system \cite{Gel:1989, Gel:1990, Gel:2008} provides, among many things, a useful tool in computing the Picard--Fuchs system for families of projective varieties. In the literature, the dif\/ferential equations obtained by the GKZ system usually factor and the Picard--Fuchs system is given by the subsystem formed by a subset of these factors. It is then natural to ask what is the reason for the factorization, and what are the geometric objects that underly the extra solutions besides the period integrals, which are integrals over the cycles in the f\/ibers of the family.

\subsection{GKZ system for the Hesse pencil}

A large part of the discussions below can be extended to slightly more general families of Calabi--Yau varieties, among which the Calabi--Yau hypersurfaces in toric varieties will be of particular interest due to their appearances in mirror symmetry. For concreteness, in the present work we shall focus on the Hesse pencil of elliptic curves as an example.

The equation of the Hesse pencil $\chi\colon \mathcal{E}\rightarrow \mathcal{B}$ is given by
\begin{gather*}
\mathcal{E}\colon \ \big\{ F(\boldsymbol{x},\psi):=x^3+y^3+z^3-3\psi x y z=0\big\}\subseteq \mathbb{P}^{2}\times \mathcal{B},
\end{gather*}
here the base $\mathcal{B}$ is a copy of $\mathbb{P}^{1}$ parametrized by $\psi$.

To def\/ine period integrals, one needs to specify a local, holomorphic section of the Hodge line bundle
\begin{gather*}
\mathcal{L}=\mathcal{R}^{0}\chi_{*}\Omega^{1}_{\mathcal{X}\,|\,\mathcal{B}}\rightarrow \mathcal{B}.
\end{gather*}
Then one can integrate the corresponding family of holomorphic top forms over the locally constant sections of a~rank~$2$ local system, which is dual to $\mathcal{R}^{1}\chi_{*}\mathbb{Z}\rightarrow \mathcal{B}$, to get period integrals.

A canonical choice for the local section is given by
\begin{gather}\label{eqnsectionHodgelinebundle}
\Omega(\psi)=\mathrm{Res} {\psi \mu_{0}\over F(\boldsymbol{x},\psi)},\qquad \mu_{0}:=zdx\wedge dy+xdy\wedge dz+ydz\wedge dx.
\end{gather}
On each f\/iber $\mathcal{E}_{\psi}$ of the family $\chi$, the $2$-form $\mu_{0}/F$ gives a meromorphic $2$-form on the ambient space $\mathbb{P}^{2}$ with a pole of order one along $\mathcal{E}_{\psi}$ and its residue gives a holomorphic top form on this elliptic curve f\/iber $\mathcal{E}_{\psi}$. The integrals of this choice of holomorphic section, over a further choice of the locally constant sections $A,B$ of the above-mentioned rank $2$ local system, gives the period integrals
\begin{gather*}
\pi_{A}(\psi)=\int_{A}\Omega(\psi),\qquad \pi_{B}(\psi)=\int_{B}\Omega(\psi).
\end{gather*}

The Picard--Fuchs equation can be derived, for example, by computing the Gauss--Manin connection or by using the Grif\/f\/iths--Dwork method. With respect to the choice $\Omega$ given above, the dif\/ferential operator annihilating the period integrals is given by
\begin{gather}\label{eqnPF}
\mathcal{L}_{\mathrm{PF}}=(\theta_{\psi}-2)(\theta_{\psi}-1)- \psi^{3}\theta_{\psi}^2,\qquad \theta_{\psi}:=\psi {\partial\over \partial \psi}.
\end{gather}
Henceforward we shall frequently use the $\theta$-operator def\/ined as above.

The details of the derivation for the GKZ system will be useful later in this work, so we recall them here following \cite{Gel:1989}.

We f\/irst extend the family a little by rewriting the equation for (the total space) of the family as
\begin{gather*}
F(\boldsymbol{x},\boldsymbol{a}):=a_{1}x^3+a_{2}y^3+a_{3}z^3+a_{0}x y z=0.
\end{gather*}

Then, one considers the actions on the polynomial $F(\boldsymbol{x},\boldsymbol{a} )$ which belong to the diagonal scalings inside the group $\mathrm{GL}_{\boldsymbol{x}}\times \mathrm{GL}_{\boldsymbol{a}}$ and hence preserve the f\/ibration structure. Here by $\mathrm{GL}_{\boldsymbol{x}}$
we mean the af\/f\/ine transformations on $\mathbb{C}^3$ parametrized by $\boldsymbol{x}=\{x,y,z\}$ space and similarly for $\mathrm{GL}_{\boldsymbol{a}}$.
Those which f\/ixes $F$ up to an overall scaling forms a subgroup~$G$. By construction, for any element in~$G$, the scaling on $\boldsymbol{a}$ is determined by that on $\boldsymbol{x}$. We can then choose the generators of $G$ to be
\begin{gather}
(x_{1},x_{2},x_{3}; a_{1},a_{2},a_{3},a_{0})\mapsto \big(\lambda x_{1}, x_{2}, x_{3}; \lambda^{-3}a_{1}, a_{2},a_{3},\lambda^{-1}a_{0}\big), \qquad \lambda\in \mathbb{C}^{*},\nonumber\\
(x_{1},x_{2},x_{3}; a_{1},a_{2},a_{3},a_{0})\mapsto \big(x_{1}, \lambda x_{2}, x_{3}; a_{1}, \lambda^{-3} a_{2},a_{3},\lambda^{-1}a_{0}\big), \qquad \lambda\in \mathbb{C}^{*},\nonumber\\
(x_{1},x_{2},x_{3}; a_{1},a_{2},a_{3},a_{0})\mapsto \big(x_{1}, x_{2}, \lambda x_{3}; a_{1}, a_{2},\lambda^{-3} a_{3},\lambda^{-1}a_{0}\big), \qquad \lambda\in \mathbb{C}^{*},\nonumber\\
(x_{1},x_{2},x_{3}; a_{1},a_{2},a_{3},a_{0})\mapsto \big( x_{1}, x_{2}, x_{3}; \lambda a_{1}, \lambda a_{2}, \lambda a_{3}, \lambda a_{0}\big), \qquad \lambda\in \mathbb{C}^{*}.\label{eqnscaling}
\end{gather}
The former three also scale the meromorphic $2$-form $\mu_{0}$ by $\lambda$, while the latter acts trivially. Hence when acting on the $2$-form
\begin{gather}\label{eqnmero2form}
\omega(\boldsymbol{a}):={\mu_{0}\over F(\boldsymbol{x},\boldsymbol{a})},
\end{gather}
the inf\/initesimal versions of these group transformations give rise to the following annihilating dif\/ferential operators called Euler homogeneity operators,
\begin{gather*}
\mathcal{Z}_{i}=\theta_{x_{i}} -3\theta_{a_{i}}- \theta_{a_{0}}-\deg_{x_{i}} \mu_{0},\qquad i=1,2,3,\nonumber\\
\mathcal{Z}_{0}=\sum_{i=1}^{3}\theta_{a_{i}}+ \theta_{a_{0}}-(-1).
\end{gather*}
Here $\deg_{x_{i}}\mu_{0}$ stands for the weight of $\mu_{0}$ under the action $x_{i}\mapsto \lambda x_{i}$, which is one in the current case.

Now the monomials in the pencil parametrized by $a_{i}$, $i=0,1,2,3$ satisfy the relation
\begin{gather*}
x_{1}^3\cdot x_{2}^3\cdot x_{3}^3=(x_{1}x_{2}x_{3})^3.
\end{gather*}
This then gives the following dif\/ferential operator that annihilates $\omega(\boldsymbol{a})$
\begin{gather*}
\mathcal{D}_{\mathrm{GKZ}}=\prod_{i=1}^{3}\partial_{a_{i}}-\partial_{a_{0}}^3.
\end{gather*}

The projection\footnote{A more intrinsic description can be given by the D-module language.}, which we denoted by $\chi_{*}$, of the above dif\/ferential operators to be base direction yields
\begin{gather*}
(\chi_{*}\mathcal{Z}_{i})\omega=0, \qquad i=0,1,2,3,\qquad (\chi_{*}\mathcal{D}_{\mathrm{GKZ}}) \omega=0.
\end{gather*}
One then translates these dif\/ferential equations to the ones satisf\/ied by the period integrals $\pi_{\gamma}(\boldsymbol{a})=\int_{\gamma}\Omega$ with respect to the holomorphic top forms $\Omega=a_{0}\omega$
\begin{gather}
\left(3\theta_{a_{i}}+\theta_{a_{0}}+\deg_{x_{i}}-1\right)\Omega=0,\qquad i=1,2,3,\nonumber\\
\left( \sum_{i=1}^3 \theta_{a_{i}}+\theta_{a_{0}}\right)\Omega,\nonumber\\
\left( {1\over \prod a_{i}}\prod \theta_{a_{i}}-{1\over a_{0}^3}(\theta_{a_{0}}-3)(\theta_{a_{0}}-2)(\theta_{a_{0}}-1) \right)\Omega=0.\label{eqnDZ}
\end{gather}
In the present case, $\deg_{x_{i}}\mu_{0}=1$, $i=1,2,3$. The former two equations allow one to make the ansatz
\begin{gather}\label{eqnansatz}
\pi_{\gamma}(\boldsymbol{a})=\pi_{\gamma}(z), \qquad z=-{a_{1}a_{2}a_{3}\over a_{0}^3}.
\end{gather}
Now when acting on a function of $z$, the last dif\/ferential equation gives
\begin{gather}\label{eqnGKZonpi}
\mathcal{D}_{\mathrm{GKZ}}\,\pi_{\gamma}=\left(\theta_{z}^3+z (-3\theta_{z}-3)(-3\theta_{z}-2)(-3\theta_{z}-1)\right)\pi_{\gamma}=0.
\end{gather}
Specializing to $a_{1}=a_{2}=a_{3}=1$, $a_{0}=-3 \psi$, by disregarding the overall constants which are irrelevant throughout the discussions, we can see that (recall~\eqref{eqnPF})
\begin{gather}
\mathcal{D}_{\mathrm{GKZ}}=\big( \theta_{\psi}^3- \psi^{-3} (\theta_{\psi}-3)(\theta_{\psi}-2)(\theta_{\psi}-1) \big)\nonumber\\
\hphantom{\mathcal{D}_{\mathrm{GKZ}}}{} =\theta_{\psi} \circ \big( \theta_{\psi}^2- \psi^{-3} (\theta_{\psi}-2)(\theta_{\psi}-1) \big)\nonumber\\
\hphantom{\mathcal{D}_{\mathrm{GKZ}}}{}=\theta_{\psi}\circ \psi^{-3}\circ \mathcal{L}_{\mathrm{PF}}.\label{eqnrightfactorization}
\end{gather}
When writing the Picard--Fuchs operator in terms of the $\alpha=\psi^{-3}$ coordinate, we use the following normalization of the leading coef\/f\/icient
\begin{gather}\label{eqnPFinalphacoordinate}
 \tilde{\mathcal{L}}_{\mathrm{PF}}= \left( \theta_{\alpha}^2- \alpha \left(\theta_{\alpha}+{1\over 3}\right)\left(\theta_{\alpha}+{2\over 3}\right) \right),
\end{gather}
so that up to a constant multiple we have $ \tilde{\mathcal{D}}_{\mathrm{GKZ}}=\theta_{\alpha}\circ \tilde{\mathcal{L}}_{\mathrm{PF}}$.

We remark that for families of hypersurface Calabi--Yau varieties in toric varieties in any dimension, similar discussions apply. In particular,
one always obtains $\mathcal{L}_{\mathrm{PF}}$ from the factor in the rightmost as in \eqref{eqnrightfactorization}.

\subsection{Calabi--Yau condition and factorization of dif\/ferential operator}

By examining the derivation, a few observations are in order. First, the f\/irst equation in \eqref{eqnDZ} is consistent with the second if and only if the following condition holds
\begin{gather}\label{eqnCYcondition}
\sum_{i=1}^{3} \deg_{x_{i}}\mu_{0}=\deg F.
\end{gather}
We call this the Calabi--Yau condition since the meromorphic form $\mu_{0}$ is a section of $\mathcal{O}_{\mathbb{P}^{2}}(-(2+1))$ and the degree of the polynomial $F$ matches with the degree of~$\mu_{0}$ exactly when $F=0$ def\/ines a Calabi--Yau hypersurface in the projective space $\mathbb{P}^2$.

Now instead of making the ansatz mentioned before in \eqref{eqnansatz}, one can eliminate the dif\/ferential operators $\theta_{a_{i}},i=1,2,3$ by solving them from the $\mathcal{Z}_{i}$-operators in \eqref{eqnDZ}. Then the $\mathcal{D}$-operator in~\eqref{eqnDZ} becomes
\begin{gather*}
\left( {1\over \prod a_{i}}\prod_{i}{\theta_{a_{0}}+\deg_{x_{i}}\mu_{0}-1\over 3}-{1\over a_{0}^3}(\theta_{a_{0}}-3)(\theta_{a_{0}}-2)(\theta_{a_{0}}-1) \right).
\end{gather*}

This dif\/ferential operator factors in the desired way when the set $\{\deg_{x_{i}}\mu_{0}-1,\, i=1,2,3\}$ has a non-empty intersection with the set $\{0,1,2\}$. Again this is trivially true in the Calabi--Yau case. For accuracy, we shall call it the right factorization to indicate that the Picard--Fuchs operator is factored out from the right. This factorization is the reason that the GKZ system gives an in-homogeneous Picard--Fuchs system.

There are natural situations where the integrand is replaced by other dif\/ferential forms with dif\/ferent scaling behaviors under the action of $G$. For example, the polynomial $F$ could be replaced by a Laurent polynomial, or the integrand by the multi-Mellin transform or the M\"ahler measure. These situations occur in local Calabi--Yau mirror symmetry~\cite{Chiang:1999tz, Hori:2000kt, Mohri:2001zz, Stienstra:2005wy} and in scattering amplitudes~\cite{Bloch:2015}. For these cases, the above procedure of deriving dif\/ferential equations from GKZ symmetries still applies.

Also for other integrands, the factorization, if exists, might be dif\/ferent. Of direct relevance to the GKZ system of the Hesse pencil is the GKZ system for the mirror geometry of $K_{\mathbb{P}^{2}}$, see~\cite{Chiang:1999tz, Hosono:2004jp}. The integrand is given by
\begin{gather}\label{eqnnoncompactvolumeform}
{1\over X_{1}X_{2}\big(a_{0}+a_{1}X_{1}+a_{2}X_{2}+a_{3}X_{1}^{-1}X_{2}^{-1}\big)+uv}{dX_{1}dX_{2}\over X_{1}X_{2}}dudv,
\end{gather}
where $(X_{1}$, $X_{2})$ are coordinates on the space $(\mathbb{C}^{*})^{2}$ and $u$, $v$ are valued in~$\mathbb{C}$. It is annihilated by
\begin{eqnarray}\label{eqnleftfactorization}
\mathcal{L}_{\mathrm{CY}_{3}} = \mathcal{L}_{\mathrm{PF}}\circ \theta_{\psi} ,\qquad \psi^{-3}=-27{a_{1}a_{2}a_{2}\over a_{0}^{3}}.
\end{eqnarray}

The combinatorial data (which is conveniently encoded in the Newton polytope or toric geometry) for its Picard--Fuchs system is identical to that of the Hesse pencil, only the scaling behavior under the symmetries in~\eqref{eqnscaling} of the integrand is dif\/ferent.

\subsection{Motivation of the work}

From the factorization in \eqref{eqnrightfactorization}, one can see that besides the period integrals, the GKZ sys\-tem~$\mathcal{D}_{\mathrm{GKZ}}$ in~\eqref{eqnGKZonpi} has one more extra solution. One of the goals of the present work is to understand this extra solution.

We also aim to understand the dif\/ference and relation between the factorizations~\eqref{eqnrightfactorization},~\eqref{eqnleftfactorization}
of the operators involved in the Hesse pencil and in the mirror geometry of~$K_{\mathbb{P}^{2}}$, respectively. That there should be such a connection is predicted by the Landau--Ginzburg/Calabi--Yau correspondence~\cite{Witten:1993}. To be a little more precise, the solutions to the GKZ system were studied in~\cite{Avram:1996} (see also~\cite{Candelas:2000}) and were identif\/ied with oscillating integrals. Hence one would expect them to appear in certain form from the perspective of the elliptic curve geometry by the Landau--Ginzburg/Calabi--Yau correspondence.

In fact, the studies in \cite{Duke:2008, Eichler:1982} imply that the extra solution to the GKZ system for the Weierstrass family can be identif\/ied with certain chain integral on the elliptic curve. As will be discussed in the present work, the chain therein is closely related to the symmetries of the Weierstrass polynomial and to certain oscillating integral. Further evidences also include some recent works~\cite{Lau:2014, Shen:2016} which suggest that part of the information encoded in the Landau--Ginzburg model should be visible in the Calabi--Yau model through the symmetries of the latter.

Finding a direct relation between the oscillating integrals in the singularity theory and (integrals of) chains living on the elliptic curves
will provide a f\/irst step towards a more conceptual understanding of the LG/CY correspondence.

\subsection*{Relation to previous works}

The explicit chain integral solutions to the GKZ system for hypersurface families were studied in~\cite{Avram:1996}. More general discussions in terms of chain integrals and D-modules were provided in the beautiful works~\cite{Bloch:2014, Huang:2015, Huang:2016}. Similar examples were discussed in~\cite{ Bloch:2015} in terms of mixed Hodge structures. These works treat the extra solution to the GKZ system as a two dimensional integral living in the ambient projective space or its blow-up. One of the main dif\/ferences between the current work and the above-mentioned ones is that we give a direct realization of the extra solution in terms of chain integral living on the elliptic curve instead of in the ambient space.

The present paper also contains several observations of\/fering connections between the extra solution to the GKZ system and some geometric objects that are of interest in mirror symmetry.

A large part of the results obtained in this work have scattered in the literature but mainly at the level of sketchy justif\/ications, our new addition on this part is then to make them more clear.

\subsection*{Outline of the paper}

In Section \ref{sec3d2dintegrals} we review the known results on the realizations of the solutions to the GKZ system in terms of $3$-dimensional oscillating integrals and $2$-dimensional chain integrals. We also interpret these integrals in terms of ones living in a non-compact Calabi--Yau variety, to incorporate the GKZ symmetries.

Section \ref{sec1dintegral} discusses the realization of the solutions to the GKZ system of the Hesse pencil in terms of objects living on the elliptic curves. First we use the Wronskian method to obtain the Eichler integral formula for the solutions. Then we express them in terms of the Beltrami dif\/ferential and cycles with vanishing period integrals. We also construct chains on the elliptic curves which give rise to the extra solution besides the period integrals.

In Section \ref{secembedding} we embed both the mirror of $K_{\mathbb{P}^{2}}$ and the elliptic curves in the Hesse pencil into some compact Calabi--Yau threefold and of\/fer a connection between the Picard--Fuchs system of the former and the GKZ system of the latter.

We conclude in Section \ref{secdiscussions} with some discussions and speculations.

\section{Invariant 3d and 2d chain integrals under GKZ symmetries} \label{sec3d2dintegrals}

The GKZ symmetries are symmetries of the polynomials $F(\boldsymbol{x},\boldsymbol{a})$, not just the varieties they def\/ine. Also the symmetries are for the forms instead of cohomology classes, as opposed to the case of the Picard--Fuchs operator derived from the Gauss--Manin connection. Hence any invariant under these symmetries will provide a solution to the resulting dif\/ferential equations.

Recall that in the above when discussing the invariance of the integrals $\pi_{\gamma}$ in \eqref{eqnansatz} under the GKZ symmetries, we used the fact that the (classes of) the cycles $\gamma$ are invariant under the scalings in~\eqref{eqnscaling}. In general, chain integrals would not satisfy the dif\/ferential equations, except when they are indeed invariant under the scalings. This will be the case when they are chains cut out by coordinate planes. This again opens the possibility that certain chain integrals could solve the GKZ system and provide extra solutions other than the cycle integrals, namely the period integrals.

\subsection{Invariant chain integrals as solutions to GKZ system}

Now we consider the so-called V-chain, see \cite{Avram:1996} and references therein, given by
\begin{gather}\label{eqn3dcontour}
D_{3}=\big\{(x,y,z)\in \mathbb{C}^3\,|\, x,y,z\geq 0\big\}\cong\mathbb{R}_{\geq 0}^3.
\end{gather}
It is indeed invariant under the transformations in~\eqref{eqnscaling}. Here we have used the coordinates $x$, $y$, $z$ in place of $x_{1}$, $x_{2}$, $x_{3}$, as we shall occasionally do throughout the work.

By applying a coordinate change, we can arrange such that $a_{i}=1$, $i=1,2,3$, $a_{0}=-3\psi$. Now we assume that the following condition
\begin{gather}\label{eqnstability}
\Re \psi\leq 0,\qquad \text{such that}\quad \Re F(x,y,z;\psi)\geq 0\quad \text{on}\quad D_{3}.
\end{gather}
The meaning of this condition will be discussed later in Remark~\ref{remsteepest}. Hence the convergence of the integral on $D_{3}$ is ensured. We can then apply the absolute convergence theorem and write
\begin{gather}
I(\psi): = \int_{D_{3}}e^{-F}\psi dxdydz\nonumber\\
\hphantom{I(\psi)}{} =\psi\sum_{k=0}^{\infty}\int_{D_{3}}
e^{-x^3}e^{-y^3}e^{-z^3}{ (3\psi)^{n}\over n!}x^n y^n z^n dxdydz
=\psi\sum_{n=0}^{\infty} { (3\psi)^{n}\over n!}{1\over3^3}\Gamma\left({n+1\over 3}\right)^3.\label{eqnoscillating}
\end{gather}
According to the residue of $n$ modulo $3$, the integral is the sum of three series
\begin{gather}\label{eq3dchainintegral}
\int_{D_{3}}e^{-F}\psi dxdydz=\sum_{i=0}^{2}\psi\sum_{k=0}^{\infty} { (3\psi)^{3k+i}\over (3k+i)!}{1\over3^3}\Gamma\left({3k+i+1\over 3}\right)^3.
\end{gather}
It breaks into the following three pieces
\begin{gather}
\psi\sum_{k=0}^{\infty} { (3\psi)^{3k}\over (3k)!}{1\over3^3}\Gamma\left({3k+1\over 3}\right)^3
={1\over 3^{3}}\cdot \left(2\pi \cdot 3^{-{1\over 2}}{\Gamma({1\over 3})^2 \over \Gamma({2\over 3})}\right)\psi\,{}_{2}F_{1}\left({1\over 3},{1\over 3};{2\over 3},\psi^{3}\right),\nonumber\\
\psi\sum_{k=0}^{\infty} { (3\psi)^{3k+1}\over (3k+1)!}{1\over3^3}\Gamma\left({3k+2\over 3}\right)^3
={1\over 3^{3}}\cdot \left(2\pi \cdot 3^{-{1\over 2}}{\Gamma({2\over 3})^2 \over \Gamma({4\over 3})}\right)\psi^2 \,{}_{2}F_{1}\left({2\over 3},{2\over 3};{4\over 3},\psi^{3}\right),\nonumber\\
\psi\sum_{k=0}^{\infty} { (3\psi)^{3k+2}\over (3k+2)!}{1\over3^3}\Gamma\left({3k+3\over 3}\right)^3
={1\over 3^{3}}\cdot \left(2\pi \cdot 3^{-{1\over 2}}{\Gamma(1)^3 \over \Gamma({4\over 3})
\Gamma({5\over 3})}\right)\psi^3 \,{}_{3}F_{2}\left(1,1,1;{4\over 3},{5\over 3};\psi^{3}\right)\!.\!\!\!\!\!\label{eqnsolsbymonodromy}
\end{gather}

There are other choices for the $V$-chain. In order for the condition $\Re F>0$ to hold and the coef\/f\/icients of~$x_{i}^3$, $i=1,2,3$ to remain, one is led to the following three chains,
\begin{gather*}
D_{3}:=C_{x}\times C_{y}\times C_{z}=(0,\infty)\times (0,\infty) \times (0,\infty),\\
\rho D_{3}:=C_{x}\times \rho C_{y}\times C_{z}=(0,\infty)\times (0,\rho\infty) \times (0,\infty),\\
\rho^2 D_{3}:=C_{x}\times \rho^2C_{y}\times C_{z}=(0,\infty)\times \big(0,\rho^2\infty\big) \times (0,\infty).
\end{gather*}
Here $\rho=\exp({2\pi i\over 3})$. Then the condition in \eqref{eqnstability} becomes
\begin{gather}\label{eqnstabilityfull}
\Re (\rho^{k}\psi)\leq 0,\qquad \text{such that}\quad \Re F(x,y,z;\psi)\geq 0\quad \text{on}\quad \rho^{k}D_{3}, \quad k=0,1,2.
\end{gather}

\begin{Remark}[steepest descent contours]\label{remsteepest}
We now make a pause and explain the condition in \eqref{eqnstabilityfull}. Note that the condition $\Re \psi\leq 0$ is not necessary in order for the chain integral $\int_{D_{3}} e^{-F}dxdydz$ to be the convergent, nor for the elliptic curve $\mathcal{E}_{\psi}=\{F(x,y,z;\psi)=0\}$ to have empty intersection with~$D_{3}$. For the former, due to the coef\/f\/icients of $x^3$, $y^3$, $z^3$, the integral is always convergent. For the latter, suppose $\mathcal{E}_{\psi}\cap D_{3}\neq \varnothing$, then one has~$\psi\in\mathbb{R}$. It is then easy to see that this is true if only and if $\psi\geq 1$ which is dif\/ferent from the condition in~\eqref{eqnstabilityfull} as well.

Hence the above condition in \eqref{eqnstabilityfull} implies the convergence condition but is stronger. In fact, any of the integrals obtained by $\rho^{k}D_{3}$ are well-def\/ined for any phase of~$\psi$ when~$\psi$ is close to~$0$, as can be seen from the explicit hypergeometric series expressions above.

For a qualitative analysis it is enough to focus on the $y$-integral part since the ranges for~$C_{x}$,~$C_{z}$ are f\/ixed and are given by the positive real axis. Hence we set $x=z=1$. Then it amounts to study the following type of Airy integral which occur in the study of the $A_{2}$-singularity theory,
\begin{gather*}
\int_{C_{\theta}} e^{-y^{3}+3\psi y }dy, \qquad C_{\theta}=e^{i\theta}\mathbb{R}_{\geq 0 }.
\end{gather*}

As already mentioned before, the process of deriving equations from symmetries can be applied to this case. In particular, the chain integrals over
$\rho^{k}C_{y}$, $k=0,1,2$ are called the so-called Scorer functions, see~\cite{Olver:2010nist}, satisfying certain 3rd order ODE.
Now for the integral to be convergent, the ray needs to sit inside one of the wedges
\begin{gather*}
W_{k}\colon \ -{\pi\over 6}+k{2\pi \over 3}\leq \arg y\leq {\pi\over 6}+k{2\pi \over 3},\qquad k=0,1,2.
\end{gather*}

Small deformations within the wedges do not af\/fect the integral, since the dif\/ference would be the integral over an arc with radius~$R$ which tends to zero as $R\rightarrow \infty$. It is in general not easy to compute the resulting integral once the chain moves out of the wedges. We hence restrict ourselves to rays inside the wedges. For the purpose of analyzing the asymptotic behavior of the integral as $\psi\rightarrow \infty$, one deforms the ray into a steepest descent contour. Among the steepest descent contours of particular importance are the ones passing through the critical points. The asymptotic expansion of such a contour integral is then completely determined from a small neighborhood of the critical point.\footnote{The asymptotic expansion derived from steepest descent method naturally leads to the so-called canonical coordinate (critical value) in singularity theory. While performing a change of variable $w=w(y;\psi)$, $y=y(w;\phi)$ such that $-y(w)^2+3\psi w= -w^3$ leads to the f\/lat coordinate $\phi$.} The steepest decent contours passing through the critical points that $C_{\theta}$ can deform to depends on the phase of $\psi$, resulting in the Stokes phenomenon.

The picture of moving integral contours to determine the asymptotics in Airy integrals also holds for the Hesse pencil case. The condition~\eqref{eqnstabilityfull} then indicates the steepest descent contours that the integral contour in consideration can deform to for the given range of~$\psi$. There are subtleties however. For example, the singularities of the GKZ system for the Hesse pencil are all regular and there are additional singularities at $\psi^3=1$.
\end{Remark}

\subsubsection{Monodromy action and functional relations}

One can also rotate the $x$, $z$ directions by powers of~$\rho$, but the resulting chains are essentially equivalent to the aforementioned three by using the actions in~\eqref{eqnscaling}. For example, the chain $\rho^{i} C_{x}\times \rho^{j} C_{y}\times C_{z}$ is equivalent to $C_{x}\times \rho^{i+j} C_{y}\times C_{z}$ as far as the integrals are concerned. These relations are nothing but a manifestation of the invariance of the integral under the action
\begin{gather}\label{eqngaugedsym}
(x,y,z;\psi)\mapsto \big(x,\lambda^{-1} y, z;\lambda \psi\big).
\end{gather}
But now due to the ``gauge f\/ixing condition" $a_{i}=1$, $i=1,2,3$, the values that $\lambda$ can take reduce from $\mathbb{C}^{*}$ to the multiplicative cyclic group $\boldsymbol{\mu}_{3}$. It is easy to see that in order to f\/ix $\psi$, the transformation must be of the form\footnote{From the perspective of the LG/CY correspondence, these symmetries should be thought of as the symmetries of the underlying Landau--Ginzburg model def\/ined at the orbifold singularity in the family~\cite{Witten:1993}.}
\begin{gather}\label{eqngaugedsymHesse}
(x,y,z)\mapsto \big(\rho^{i}x,\rho^{j}y,\rho^{k}z\big),\qquad i+j+k\equiv 0~\mod~3.
\end{gather}

According to the invariance under \eqref{eqngaugedsym}, the integrals over $\rho^{k}D_{3}$ then satisfy the functional relations
\begin{gather}
I_{\rho}(\psi): = \int_{\rho D_{3}}e^{-F}\psi dxdydz= I(\rho \psi)=\rho J_{1}(\psi)+\rho^2 J_{2}(\psi)+J_{3}(\psi),\nonumber\\
I_{\rho^2}(\psi): = \int_{\rho^2 D_{3}}e^{-F}\psi dxdydz= I\big(\rho^2 \psi\big)=\rho^2 J_{1}(\psi)+\rho J_{2}(\psi)+J_{3}(\psi).\label{eqnfunctionaleqn}
\end{gather}

On the other hand, the solutions annihilated by $\mathcal{L}_{\mathrm{GKZ}}$ in~\eqref{eqnrightfactorization} are easily seen to be
\begin{gather}
\psi \,{}_{3}F_{2}\left({1\over 3},{1\over 3},{1\over 3};{1\over 3},{2\over 3};\psi^{3}\right)=\psi \,{}_{2}F_{1}\left({1\over 3},{1\over 3};{2\over 3};\psi^{3}\right),\nonumber\\
 \psi^2 \,{}_{3}F_{2}\left({2\over 3},{2\over 3},{2\over 3};{2\over 3},{4\over 3};\psi^{3}\right)=\psi^2 \,{}_{2}F_{1}\left({2\over 3},{2\over 3};{4\over 3};\psi^{3}\right),\nonumber\\
\psi^3 \,{}_{3}F_{2}\left(1,1,1;{4\over 3},{5\over 3};\psi^{3}\right).\label{eqnexplicitGKZhypersols}
\end{gather}
By comparing these with the above three chain integrals $I(\psi)$, $I_{\rho}(\psi)$, $I_{\rho^2}(\psi)$, we can see indeed the 3d chain integrals give the full set of solutions to the GKZ system.

\subsubsection{Period integrals as dif\/ferences of chain integrals}\label{secdifferencesofchainintegrals}

Recall from \eqref{eqnrightfactorization} that the period integrals are solutions annihilated by~$\mathcal{L}_{\mathrm{PF}}$ and hence are given by
\begin{gather}\label{eqnperiods}
\pi_{1}(\psi)=\psi \,{}_{2}F_{1}\left({1\over 3},{1\over 3};{2\over 3};\psi^{3}\right),\qquad
\pi_{2}(\psi)=\psi^2 \,{}_{2}F_{1}\left({2\over 3},{2\over 3};{4\over 3};\psi^{3}\right).
\end{gather}
These are proportional to the solutions to the GKZ system which are given in the f\/irst two in~\eqref{eqnsolsbymonodromy}. One can also check directly that for the extra solution to $\mathcal{L}_{\mathrm{GKZ}}$ in~\eqref{eqnexplicitGKZhypersols} one has
\begin{gather*}
\psi^{-3}\mathcal{L}_{\mathrm{PF}} \left(\psi^3 \,{} _{3}F_{2}\left(1,1,1;{4\over 3},{5\over 3};\psi^{3}\right)\right)={2\over 9}.
\end{gather*}

\begin{Remark}
A more convenient choice of basis (for the integrality of the connection matrices) near this point is given by~\cite{Erdelyi:1981}
\begin{gather*}
\tilde{\pi}_{1}=-\rho {\Gamma({1\over 3}) \over \Gamma({2\over 3})^2}
\psi \,{}_{2}F_{1}\left({1\over 3},{1\over 3};{2\over 3};\psi^{3}\right), \qquad
\tilde{\pi}_{2}=\rho^2 {\Gamma(-{1\over 3}) \over \Gamma({1\over 3})^2} \psi^2 \,{}_{2}F_{1}\left({2\over 3},{2\over 3};{4\over 3};\psi^{3}\right).
\end{gather*}
In terms of the parameter $\alpha=\psi^{-3}$, the singularities of the Hesse pencil include the cusp singularities $\alpha=0,1$ and the orbifold singularity $\alpha=\infty$. The periods corresponding to the vanishing cycles at $\psi^{-3}=0$, $\psi^{-3}=1$ are given by,
\begin{gather}
\omega_{0}= {}_{2}F_{1}\left({1\over 3},{2\over 3};1;\psi^{-3}\right)=\tilde{\pi}_{1}+\tilde{\pi}_{2},\nonumber\\
\omega_{1}={i\over \sqrt{3}}\,{}_{2}F_{1}\left({1\over 3},{2\over 3};1;1-\psi^{-3}\right)={i\over \sqrt{3}}\big({-}\rho \tilde{\pi}_{1}+\rho^2\tilde{\pi}_{2}\big).\label{eqnanalyticcontinuation}
\end{gather}
See \cite{Shen:2016} for a collection of results.
\end{Remark}

The solutions to the GKZ system are naturally expanded around the orbifold point $\psi=0$ in the base $\mathcal{B}$. If we look at the monodromy around this orbifold point, the period integrals (i.e., solutions to the Picard--Fuchs system) correspond to the f\/irst two in~\eqref{eqnsolsbymonodromy} which have non-trivial monodromies under the action $\psi\mapsto e^{2\pi i}\psi$. The extra solution~$J_{3}$ is the monodromy invariant one which is therefore invisible from the Picard--Fuchs equation. Recall that the local monodromy action near the orbifold point is rooted in the ``gauged symmetry''~in~\eqref{eqngaugedsymHesse} and is what leads to the functional relations in~\eqref{eqnfunctionaleqn}. All these suggest that the orbifold singularity plays a special role and can detect more information about the family other than the vanishing cycles which are topological. We shall say more about this later in Section~\ref{sec1dintegral}.

\begin{figure}[h]\centering
\includegraphics[scale=0.5]{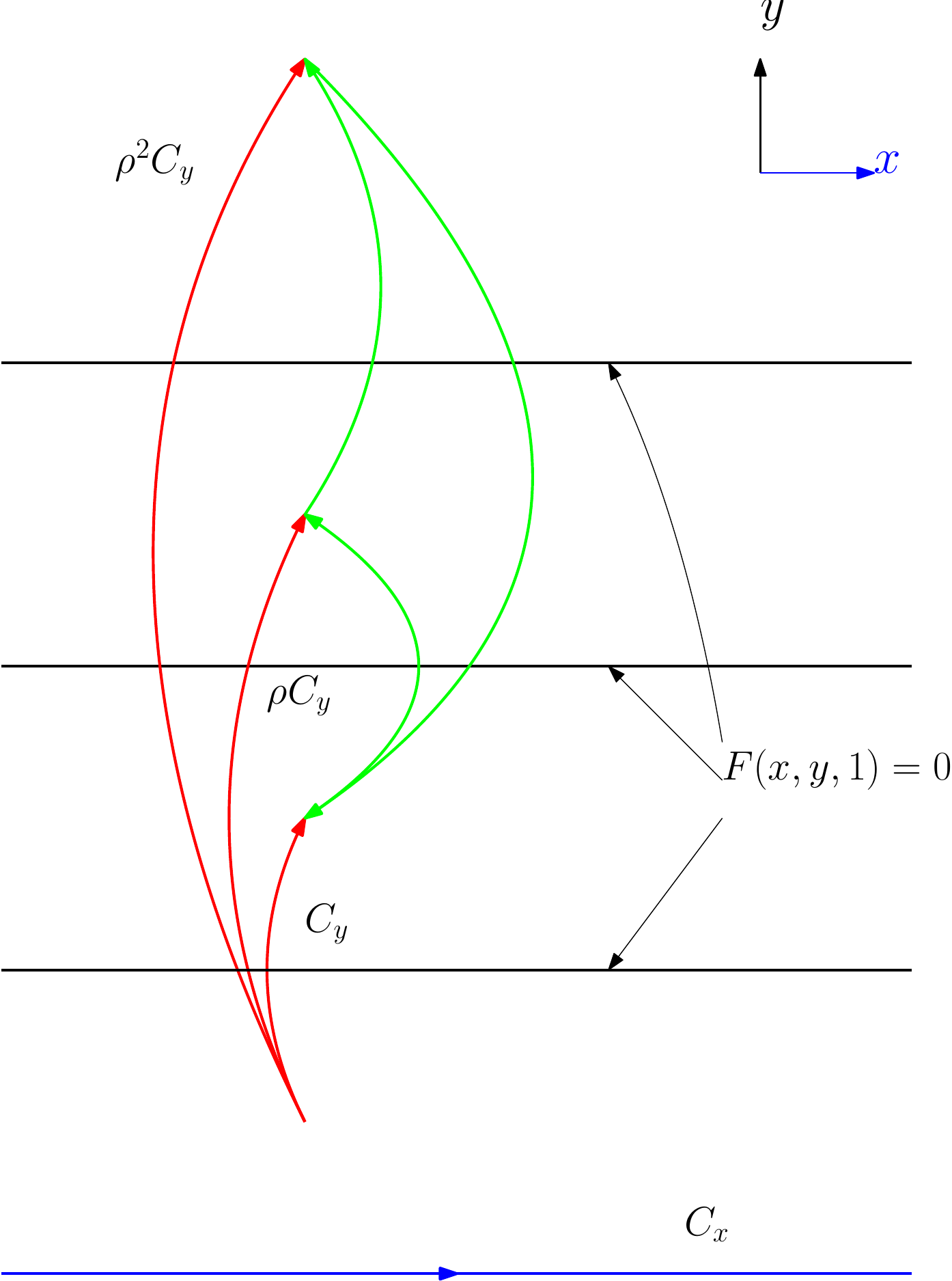}
 \caption{Dif\/ferences of 2-chains give tubular neighorhoods.} \label{figuredifferenceofchains}
\end{figure}

The dif\/ferences between the 3d chain integrals give rise to cycle integrals
\begin{gather*}
I_{\rho}(\psi)-I(\psi) = \int_{\rho D_{3}-D_{3}}e^{-F}dxdydz =(\rho-1)J_{1}+\big(\rho^2-1\big)J_{2},\\
I_{\rho^2}(\psi)-I_{\rho}(\psi) = \int_{\rho^2 D_{3}-\rho D_{3}}e^{-F}dxdydz =\big(\rho^2-\rho\big)J_{1}+\big(\rho-\rho^2\big)J_{2}.
\end{gather*}
One can also check directly that these cycles correspond to cycles on the elliptic curves without using the relations to the periods in~\eqref{eqnperiods}. To do this we note that in the common region of~$\psi$ such that for both~$I_{\rho^{k_{1}}}(\psi)$ and~$I_{\rho^{k_{2}}}(\psi)$, $k_{1}\neq k_{2}~\mathrm{mod}~3$ the condition~\eqref{eqnstabilityfull} holds, the chains $C_{x}\times \rho^{k_{1}} C_{y}$, $C_{x}\times \rho^{k_{2}} C_{y}$ have no intersection with the elliptic curve def\/ined by~$F=0$. The dif\/ference gives a tubular neighborhood of a certain branch $C_{k_{1},k_{2}}$ of $\{(x,y,z)\in \{F=0\}\,|\, x\in C_{x}\}$. Then using the residue calculus, one f\/inds a chain integral on the elliptic curve
\begin{gather}\label{eqn1dintegral}
I_{\rho^{k_{1}}}(\psi)-I_{\rho^{k_{2}}}(\psi)=\int_{C_{x}\times (\rho^{k_{1}} C_{y}-\rho^{k_{2}} C_{y})} {\psi \mu_{0}\over F}
=\int_{C_{x}}\mathrm{Res} {\psi \mu_{0}\over F}\Big|_{C_{k_{1},k_{2}}},\qquad k=0,1,2.
\end{gather}
See Fig.~\ref{figuredifferenceofchains} for an illustration.

It is a classical result, see \cite{Artebani:2006,Dolgachev:1997}, that the Hesse pencil arise as the equivariant embedding of elliptic curves with $3$-torsion structure to the projective plane via the theta functions. Here equivariance means that the action of translations by the group $E[3]$ of $3$-torsion points (which are lattice points in ${1\over 3} (\mathbb{Z}\oplus \mathbb{Z}\tau)$ if the elliptic curve is realized as $E=\mathbb{C}/(\mathbb{Z}\oplus \mathbb{Z}\tau)$) on the curve~$E$ gets mapped to projective transformations on~$\mathbb{P}^{2}$. Using the particular choices for the theta functions in~\cite{Dolgachev:1997}, these projective transformations are
\begin{gather}
\sigma_{1}={1\over 3}\colon \quad (x,y,z)\mapsto \big(x,\rho y,\rho^2 z\big),\nonumber\\
\sigma_{2}={\tau\over 3}\colon \quad (x,y,z)\mapsto (y,z,x).\label{eqntranslationaction}
\end{gather}

The zeros of the coordinate functions $x$, $y$, $z$ correspond to those of the theta functions used to def\/ine the embedding. Hence the dif\/ference between any two of them (in particular the endpoints of $C_{k_{1},k_{2}}$) are nothing but the $3$-torsion points on the elliptic curve. Using the translations in~\eqref{eqntranslationaction} which f\/ixes $\psi$, the chain $C_{k_{1},k_{2}}$ in~\eqref{eqn1dintegral} can then produce full cycles. This shows that the dif\/ference of the chain integrals given in~\eqref{eqn1dintegral} essentially give the period integrals.

\subsection{Integrals on a local Calabi--Yau}\label{seclocalCYintegrals}

In the above we explained that the 3d chain integrals in~\eqref{eq3dchainintegral} give the full set of solutions to the GKZ system. These integrals are the so-called oscillating integrals on~$\mathbb{C}^3$
\begin{gather}\label{eqnosscilatingintegralaffine}
\int a_{0} e^{-F}dxdydz.
\end{gather}
Here we have omitted the chains in the integration since by construction they are invariant under the GKZ symmetries in~\eqref{eqnscaling} and
are not important in the discussions below.

The integrand is invariant under the symmetries in \eqref{eqnscaling} which f\/ix the polynomial $F(\boldsymbol{x},\boldsymbol{a})$. However, they are not invariant under the symmetries which do not f\/ix~$F$ but result in scalings on~$F$. This set of symmetries is generated by
\begin{gather}\label{eqnanomaloustrans1}
(\boldsymbol{x},\boldsymbol{a})\mapsto (\boldsymbol{x},\lambda \boldsymbol{a}),
\end{gather}
and
\begin{gather}\label{eqnanomaloustrans2}
(\boldsymbol{x},\boldsymbol{a})\mapsto (\lambda \boldsymbol{x},\boldsymbol{a}).
\end{gather}
For example, the former gives the equation
\begin{gather*}
\sum_{i=0}^{3} \theta_{a_{i}} F=F,\qquad \sum_{i=0}^{3} \theta_{a_{i}} (a_{0}dxdydz)=(a_{0}dxdydz).
\end{gather*}
The corresponding operators act as the Euler operators on homogeneous functions. Unlike the $\mu_{0}/F$ case in~\eqref{eqnmero2form}, one does not have a symmetry on the integrand since
\begin{gather*}
\sum_{i=0}^{3} \theta_{a_{i}} \big(e^{-F}a_{0} dxdydz \big)= - \big(e^{-F} F a_{0}dxdydz\big)\neq - \big(e^{-F} a_{0}dxdydz\big).
\end{gather*}
That is, the transformations in \eqref{eqnanomaloustrans1}, \eqref{eqnanomaloustrans2}, which are redundant for $\mu_{0}/F$ when the Calabi--Yau condition~\eqref{eqnCYcondition} holds, do not seem to yield symmetries for~$e^{-F}dxdydz$. However, one knows from the explicit computations that the chain integrals in \eqref{eqnosscilatingintegralaffine} do generate the full space of solutions to the GKZ system and hence should be invariant under these transformations.

To resolve this conf\/lict, we are led to the following more correct interpretation of the oscillating integral in~\eqref{eqnosscilatingintegralaffine}. First we note that the above two transformations \eqref{eqnanomaloustrans1}, \eqref{eqnanomaloustrans2} are related by a transformation which does preserve $F$, hence we only need to consider one of them. For simplicity, we focus on the latter.
Motivated by~\cite{Witten:1993}, we think of the above integral as one on the total space of~$K_{\mathbb{P}^{2}}$. We choose~$s(\boldsymbol{a})$ to be coordinate on the f\/iber with respect to the trivialization~$a_{0}\mu_{0}$. Note that $a_{0}\mu_{0}$ fails to represent a nonzero section precisely at the orbifold point~$a_{0}=0$ on the base of the elliptic curve family. This is why the coordinate~$s(\boldsymbol{a})$ is moduli dependent in order to render $s(\boldsymbol{a}) a_{0}\mu_{0}$ well-def\/ined. Then the holomorphic top form on $K_{\mathbb{P}^{2}}$ is $a_{0}\mu_{0}\wedge ds(\boldsymbol{a})$. Now we consider the dif\/ferential form on the CY threefold $K_{\mathbb{P}^{2}}$
\begin{gather*}
e^{-sF(\boldsymbol{x},\boldsymbol{a})}a_{0}\mu_{0}\wedge ds(\boldsymbol{a}) .
\end{gather*}
Since $s(\boldsymbol{a})a_{0}\mu_{0}$ gives a meromorphic section of $K_{\mathbb{P}^{2}}\rightarrow \mathbb{P}^2$, under the $\mathbb{C}^{*}$-actions in~\eqref{eqnanomaloustrans1},~\eqref{eqnanomaloustrans2} the quantity $s(\boldsymbol{a}) F$ is invariant.

We regard $W=sF$ as a function in the coordinate ring of the variety $K_{\mathbb{P}^{2}}$. It follows that the Calabi--Yau condition~\eqref{eqnCYcondition} simply means that $W$ is homogeneous of degree one in the f\/iber coordinate~$s$
\begin{gather*}
\nu:=\deg_{s}(W)=1.
\end{gather*}
This way of looking at the scaling behavior of the form $\mu_{0}$ is convenient, especially when there is no term $a_{0}\prod x_{i}$ involved in the polynomial $F$ which was used to absorb the shift in \eqref{eqnDZ} that comes from the action on the $\mu_{0}$ part.

One can then compute the resulting integral as follows
\begin{gather}\label{eqndimreduction}
\int_{0}^{\infty} \int e^{-sF}a_{0}\mu_{0}\wedge ds= \int {1\over F}a_{0}\mu_{0}.
\end{gather}
In particular, in the patch $z=1$, one has
 \begin{gather*}
\int_{0}^{\infty} \int e^{-(sz^3) {F\over z^3}}a_{0} d\left({x\over z}\right)\wedge d\left({y\over z}\right)\wedge ds.
\end{gather*}
Now one can formally make the change of variable $sz^3\mapsto z^3$, as a computational shortcut, then the above integral becomes
 \begin{gather*}
\int_{0}^{\infty} \int e^{-F}3a_{0}\, dx\wedge dy\wedge dz.
\end{gather*}
This gives the oscillating integrals discussed earlier in~\eqref{eqnosscilatingintegralaffine}.

In sum, in order to respect all of the GKZ symmetries, the integral in \eqref{eqnosscilatingintegralaffine} should be interpreted as one on~$K_{\mathbb{P}^{2}}$. In doing actual computations, we shall however think of the integral as if it is on $\mathbb{C}^3$ for convenience.

For the 3d chain integral in \eqref{eq3dchainintegral}, according to \eqref{eqndimreduction}, one gets the 2d real integral which in the af\/f\/ine coordinate $z=1$ becomes
\begin{gather}\label{eqn2drealintegral}
\int_{0}^{\infty}\int_{0}^{\infty} {a_{0}dxdy\over F(x,y,1;\psi )}.
\end{gather}
The resulting 2d chain is interpreted in \cite{Huang:2015} as an element in the relative homology $H_{2}(\mathbb{P}^{2}-E$, $\Delta-E\cap \Delta)$, where $\Delta=\{xyz=0\}$. Similar examples are discussed in~\cite{Bloch:2015} in which the pencil of cubic curves have base points lying on the integral domain and a blow-up is needed.

\section{Chains on the elliptic curves and orbifold singularities}\label{sec1dintegral}

It will be more satisfactory if one can f\/ind chains in the elliptic curve f\/ibers that give rise to the extra solutions to the GKZ system. As mentioned above, this will then establish a link between the oscillating integrals~\eqref{eqnosscilatingintegralaffine} in the singularity theory and objects in the elliptic curve geometry. Since the integral contour in~\eqref{eqn2drealintegral} is not a tubular neighborhood of a chain on the elliptic curve, a direct dimension reduction is not available.

Instead, we shall f\/irst derive an integral formula for the extra solution basing on the in-homogeneous Picard--Fuchs equation and the Wronskian method. The relation to modular forms, which is special in the current example, gives an Eichler integral. Also the integral formula of\/fers a nice interpretation of the extra solution in terms of the Beltrami dif\/ferential which captures the deformation of complex structures.

Independently, we obtain a chain integral on the elliptic curve for the extra solution, motivated by the special role played by the orbifold singularity in the moduli space.

\subsection{Wronskian method: Eichler integral}\label{secWronskian}

We use the Wronskian method to obtain an Eichler integral formula for the solution $I(\psi)$ following~\cite{Bloch:2015, Duke:2008, Eichler:1982}. Recall that the extra solution $I(\psi)$ to $\mathcal{D}_{\mathrm{GKZ}}$ in~\eqref{eqnrightfactorization} must solve the in-homogeneous Picard--Fuchs equation
\begin{gather*}
\big( \theta_{\psi}^2- \psi^{-3} (\theta_{\psi}-2)(\theta_{\psi}-1) \big)=C
\end{gather*}
for some constant $C$. Taking any basis of the periods $u_{1}$, $u_{2}$ annihilated by the Picard--Fuchs operator $\mathcal{L}_{\mathrm{PF}}$ in~\eqref{eqnPF}, then according to the standard Wronskian method one has
\begin{Theorem} The solutions $I(\psi)$ to the GKZ system for the Hesse pencil are given by
\begin{gather}\label{eqnWronskian}
I(\psi)=a u_{1}(\psi)+b u_{2}(\psi)+c \int^{\psi} {1\over (1-v^{-3})v^{2}}{1\over W(v)} (u_{1}(\psi) u_{2}(v)-u_{2}(\psi) u_{1}(v))dv,
\end{gather}
for some constants $a$, $b$, $c $.
\end{Theorem}
The lower bound in the integral does not matter: two dif\/ferent choices for the lower bound result in a change on $a$, $b$.

The Wronskian
\begin{gather*}
W(\psi)= (u_{1}'(\psi)u_{2}(\psi)-u_{1}(\psi)u_{2}'(\psi))
\end{gather*}
can be easily computed by using the Schwarzian of the Picard--Fuchs equation.

It is known that the Hesse pencil is parametrized by the modular curve $\Gamma_{0}(3)\backslash \mathcal{H}^{*}$ whose Hauptmodul $\alpha(\tau)$ can be found, e.g., in~\cite{Maier:2009}. We take the basis $u_{1}$, $u_{2}$ to be the pe\-riods $ \omega_{0}(\alpha)$,~$\omega_{1}(\alpha)$ near the inf\/inity cusp given in~\eqref{eqnanalyticcontinuation}, with~\cite{Berndt:1995} $\tau=\omega_{1}/\omega_{0}$.
See~\cite{Shen:2016} for a collection of the formulas. Now the last term in~\eqref{eqnWronskian} is
\begin{gather*}
\int^{\alpha}{1\over v^2(1-v)}{1\over W(v)} ( \omega_{0}(v)\omega_{1}(\alpha)-\omega_{0}(\alpha)\omega_{1}(v))dv\\
\qquad {} =\omega_{0}(\alpha)\int^{\alpha} {1\over v^2(1-v)}{1\over W(v)} \omega_{0}(v)(\tau(\alpha)-\tau(v))dv,
\end{gather*}
up to a constant multiple. Then we get
\begin{Corollary}
Denote the normalized period $I/\omega_{0}$ by $t_{\mathrm{GKZ}}$, then one has the following Eichler integral expression of~$t_{\mathrm{GKZ}}$ near the infinity cusp
\begin{gather}
t_{\mathrm{GKZ}} = a+b\tau+c\int^{\alpha} {1\over v^2(1-v)}{1\over W(v)}\omega_{0}(v)(\tau(\alpha)-\tau(v))dv\,\nonumber\\
\hphantom{t_{\mathrm{GKZ}}}{} = a+b\tau+c \int^{\tau} (1-\alpha(v))\omega_{0}^3(v)(\tau-v)dv,\label{eqnnormalizedperiod}
\end{gather}
for some constants $a$, $b$, $c$.
\end{Corollary}
The above formula in \eqref{eqnnormalizedperiod} is consistent with the result that
\begin{gather*}
\tilde{\mathcal{L}}_{\mathrm{GKZ}}(\omega_{0}t_{\mathrm{GKZ}})
=\theta_{\alpha}\circ \tilde{\mathcal{L}}_{\mathrm{PF}}(\omega_{0}t_{\mathrm{GKZ}})
=\theta_{\alpha}\circ {1\over (1-\alpha)\omega_{0}^{3}}\circ \partial_{\tau}^{2}t_{\mathrm{GKZ}}=0.
\end{gather*}

Dif\/ferent choices for the reference point in the Eicher integral will af\/fect the last term by a~quantity whose second derivative in~$\tau$ vanishes and hence is a period integral.

We now relate the solutions to modular forms. From~\eqref{eqnnormalizedperiod} it follows that
\begin{gather}\label{eqn2ndderivativeoft}
\partial_{\tau}^2 t_{\mathrm{GKZ}}=c (1-\alpha)\omega_{0}^3,
\end{gather}
for some constant $c$. Moreover, the quantity $(1-\alpha)\omega_{0}^3$ is equal to the following modular form of weight $3$ for the modular group~$\Gamma_{0}(3)$
\begin{gather*}
B(\tau)={\eta(\tau)^3\over \eta(3\tau)},
\end{gather*}
see \cite{Maier:2009} for details. See also~\cite{Zhou:2013hpa} for detailed discussions on the computations on periods. In~\eqref{eqn2ndderivativeoft}, when $c=0$ one gets the period integrals, otherwise one gets the extra solution to the GKZ system. For simplicity, we set $c=1$ below. The modular form $B^{3}$ has a nice Eisenstein series and hence Lambert series formula given by
\begin{gather*}
B^{3}(\tau)=1-9\sum_{n\geq 1}\chi_{-3}(n) {n^2 q^{n}\over 1-q^{n}}, \qquad q=\exp (2\pi i\tau).
\end{gather*}
Here $\chi_{-3}$ is the Legendre symbol which takes the values $0$, $1$, $-1$ on integers of the form $3k$, $3k+1$, $3k+2$, respectively.
Hence we obtain
\begin{gather}\label{eqnsolintermsofq}
t_{\mathrm{GKZ}}={1\over 2}\tau^2+b\tau+a+9\sum_{n\geq 1}\chi_{-3}(n) \mathrm{Li}_{2}(q^{n}).
\end{gather}

\begin{Remark}The normalized period $t_{\mathrm{GKZ}}$ should be contrasted to the f\/lat coordinate $t$ for the mirror of the A-model geometry of $K_{\mathbb{P}^{2}}$, which is a normalized period solved from the 3rd order Picard--Fuchs equation in~\eqref{eqnleftfactorization} and arises as the integral of the Mahler measure~\cite{Stienstra:2005wy}. It satisf\/ies $\theta_{\alpha}t=\omega_{0}$. By using the Schwarzian this becomes
\begin{gather}\label{eqn1stderivativeoftCY3}
\partial_{\tau} t=c B(\tau)^{3},
\end{gather}
for some constant $c$. By using its expected boundary behavior, one gets
\begin{gather*}
e^{t}=-q\prod_{n\geq 1}(1-q^{n})^{9n\chi_{-3}(n) }, \qquad q=e^{2\pi i \tau}.
\end{gather*}
The inversion of this quantity carries interesting enumerative meaning in Gromov--Witten theory. See~\cite{Mohri:2001zz, Stienstra:2005wy, Zhou:2014thesis} for detailed discussions.

By comparing \eqref{eqn2ndderivativeoft} with \eqref{eqn1stderivativeoftCY3}, and using the properties of the special geometry~\cite{Strominger:1990pd} on the moduli space, one can see that~$t_{\mathrm{GKZ}}$ is actually related to the quantum volumes of cycles~\cite{Hosono:2004jp} in the A-model Calabi--Yau geometry under mirror symmetry. To be a little more detailed, denoting the prepotential by~$F(t)$, then the quantum volumes are given by the normalized periods $1$, $t$, $\partial_{t}F(t)$, $2F(t)-t \partial_{t}F(t)$. The normalized solutions to the GKZ system are then, up to unimportant terms,
\begin{gather*}
1=\partial_{t}(t), \qquad \tau=\partial_{t}(\partial_{t}F(t)),\\
t_{\mathrm{GKZ}}(\tau)=-\partial_{t}\left(2F(t)-t \partial_{t}F(t)\right)=t\partial_{t}(\partial_{t}F(t))-\partial_{t}F(t).
\end{gather*}
An amusing observation is that $t_{\mathrm{GKZ}}(\tau)$ is the Legendre dual of $\partial_{t}F(t)$ and vice versa.

Since in the current example the Yukawa coupling, which is def\/ined to be $\partial_{t}^3 F(t)=\partial_{t}\tau$, is non-vanishing, we can write a derivatives in $\tau$ in terms of that in $t$.
Ignoring the overall multiplicative factors, and focusing on the normalized periods, we get the simplif\/ications
\begin{gather}
\mathcal{D}_{\mathrm{GKZ}} = \partial_{\tau}\circ {\partial\tau\over \partial t}\circ \partial_{\tau}^2
=\partial_{\tau}\partial_{t}\partial_{\tau}\sim \partial_{t}\circ \partial_{t}\partial_{\tau},\\
\mathcal{L}_{\mathrm{CY}^{3}} = \partial_{t}\circ {\partial t\over \partial \tau}\circ \partial_{t} \circ {\partial}_{t}
=\partial_{t}\partial_{\tau}\circ \partial_{t}.
\end{gather}
We shall say more about the relation between them in Section~\ref{secembedding}.
\end{Remark}

\begin{Remark} Since the oscillating integral $\omega_{0} t_{\mathrm{GKZ}}$ appears naturally in the Landau--Ginzburg B-model, in particular, through the Frobenius manifold structure, it is natural to ask whether the normalized period $t_{\mathrm{GKZ}}$ is also related to the enumerative geometry of the mirror LG A-model, similar to the f\/lat coordinate $t$ for the mirror of the A-model geometry of~$K_{\mathbb{P}^{2}}$.
\end{Remark}

The solution displayed in \eqref{eqnsolintermsofq} agrees with the fact that the solutions of the GKZ system must contain a solution with $\log^{2}\alpha$ behavior. The latter ref\/lects that the indicial equation has three roots $0$, $0$, $0$ at the point $\alpha=0$ (around which a basis of solutions can be obtained via Frobenius method). The indeterminacy $a$, $b$ indicates that the extra solution is subject to addition by the other two solution which are periods and do not af\/fect the $\log^2 \alpha$ behavior. For a~given solution, say~$J_{3}$, the constants $a$, $b$, $c$ in~\eqref{eqnnormalizedperiod} can be f\/ixed following the standard method. We shall not do this here. Instead, we discuss the representation of the extra solution near the orbifold point $\psi=0$ around which qualitatively analyzing the solutions is convenient since the local monodromy action can be diagonalized.

We compute the Wronskian and get,
\begin{gather*}
W(\psi)=\psi^2 \big(1-\psi^3\big)^{-1}.
\end{gather*}
We take the basis of solutions $u_{1}$, $u_{2}$ to be the ones $\pi_{1}$, $\pi_{2}$ in~\eqref{eqnperiods}. Then we obtain
\begin{gather*}
I(\psi) =a \pi_{1}(\psi)+b \pi_{2}(\psi)+c\int^{\psi} v^{-1}( \pi_{1}(\psi) \pi_{2}(v)- \pi_{2}(\psi) \pi_{1}(v))dv,
\end{gather*}
The local uniformizing variable near the orbifold point $\psi=0$ on the base $\mathcal{B}$ can be taken to be $s=\pi_{2}/\pi_{1}$. Then in terms of $s$ one has
\begin{Corollary}
The local expansion of the solutions to the GKZ system for the Hesse pencil near the orbifold point is given by
\begin{gather*}
I(s) =a \pi_{1}(s)+b \pi_{2}(s)+c\int^{s} v^{-1}( \pi_{1}(s) v \pi_{2}(v)-s \pi_{1}(s) \pi_{1}(v)){d\psi\over ds}(v)dv\\
\hphantom{I(s)}{} = a \pi_{1}(s)+b s \pi_{1}(s) + c\pi_{1}(s)\int^{s} v^{-1} \pi_{1}(v)( v -s ){d\psi\over ds}(v)dv.
\end{gather*}
\end{Corollary}

Now it suf\/f\/ices to discuss the integral $J_{3}$ in terms of the above form since the other two solutions are period integrals which are solutions to the homogeneous Picard--Fuchs equation. Hence we want to determine the constants $a$, $b$, $c$ in the equality
\begin{gather*}
\psi^{3}\,{}_{3}F_{2}\left(1,1,1;{4\over 3},{5\over 3};\psi^3\right)=
a \psi\,{}_{2}F_{1}\left({1\over 3},{1\over 3};{2\over 3};\psi^3\right)+
b \psi^2\,{}_{2}F_{1}\left({1\over 3},{1\over 3};{2\over 3};\psi^3\right)\\
\hphantom{\psi^{3}\,{}_{3}F_{2}\left(1,1,1;{4\over 3},{5\over 3};\psi^3\right)=}{} +c \int_{0}^{\psi} v^{-1} \left(\psi v^2\,{}_{2}F_{1}\left({1\over 3},{1\over 3};{2\over 3};\psi^3\right)
\,{}_{2}F_{1}\left({2\over 3},{2\over 3};{4\over 3};v^3\right)\right.\\
\left.\hphantom{\psi^{3}\,{}_{3}F_{2}\left(1,1,1;{4\over 3},{5\over 3};\psi^3\right)=}{} -\psi^2 v\,{}_{2}F_{1}\left({1\over 3},{1\over 3};{2\over 3};v^3\right)
\,{}_{2}F_{1}\left({2\over 3},{2\over 3};{4\over 3};\psi^3\right) \right)dv.
\end{gather*}
We then use the series formula for hypergeometric functions. Without doing any calculations, we can see that due to the monodromy behavior near~$\psi=0$, we must have $a=b=0$. Then by comparing the coef\/f\/icients of~$\psi^3$, we are led to
\begin{gather*}
c= -2.
\end{gather*}
Therefore, the in-homogeneous contribution in the solution in terms of the Wronskian gives the monodromy invariant chain integral~$J_{3}$.

Again this approach singles out the special role of the orbifold point where the gauged symmetry in \eqref{eqngaugedsymHesse} results in the monodromy (under which the solutions have dif\/ferent behaviors).

\subsection{Wronskian method: vanishing periods and Beltrami dif\/ferential}

We now give a geometric interpretation of the last term in~\eqref{eqnWronskian} obtained by the Wronskian method
\begin{gather*}
J(\alpha):=\int^{\alpha}{1\over v} (u_{1}(v)u_{2}(\alpha)-u_{1}(\alpha)u_{2}(v))dv.
\end{gather*}
This naturally lives in the homology of the total space of the elliptic curve f\/ibration, similar to the integral over the Lefschetz thimble.

To see this, we choose $\gamma_{1}$, $\gamma_{2}$ to be any locally constant basis of $H_{1}(E,\mathbb{Z})$ which can be thought of as coming from the marking $m\colon H_{1}(E,\mathbb{Z})\cong \mathbb{Z}^{2}$ for a generic reference f\/iber. We again take~$\Omega(v)$ to be the section of the Hodge line bundle specif\/ied in~\eqref{eqnsectionHodgelinebundle}. The period integrals over the two cycles $\gamma_{1}$, $\gamma_{2}$, with respect to~$\Omega(v)$, gives a basis $u_{1}(v)$, $u_{2}(v)$ of solutions to the Picard--Fuchs equation. Then we rewrite $J(\alpha)$ as
\begin{gather*}
J(\alpha)=\int^{\alpha} {1\over v}dv\int_{\gamma(v;\alpha)}\Omega(v),\qquad
\gamma(v;\alpha)=u_{2}(\alpha)\gamma_{1}-u_{1}(\alpha)\gamma_{2}.
\end{gather*}

Fixing $\alpha$, the cycle $\gamma(v;\alpha)$ is locally constant in $v$ due to parallel transport. It is singled out, up to a constant multiple, by the condition
\begin{gather}\label{eqnvanishingperiod}
\int_{\gamma(\alpha;\alpha)}\Omega(\alpha)=0.
\end{gather}
That is, away from the orbifold point in the moduli space, it is exactly the unique cycle in~$H_{1}(\mathcal{E}_{\alpha},\mathbb{C})$ which is the Poincar\'e dual of $\Omega(\alpha)$ and hence gives the vanishing period in the f\/iber~$\mathcal{E}_{\alpha}$.

It is easy to check that the cycle $\gamma(v;\alpha)$ is independent of the marking and in particular is invariant under monodromy. It also varies holomorphically in $\alpha$. We call it the \emph{singular cycle}. Note the dif\/ference between the singular cycles and vanishing cycles (def\/ined with respect to the cusps).

It follows that the quantity $J(\alpha)$ measures the area of the 2-dimensional region ``swept out'' by the singular cycle at the point $\alpha$ through parallel transport, with respect to the holomorphic volume form ${dv\over v}\wedge \Omega(v)$ on the total space of the f\/ibration.

\begin{figure}[h]\centering
\includegraphics[scale=0.66]{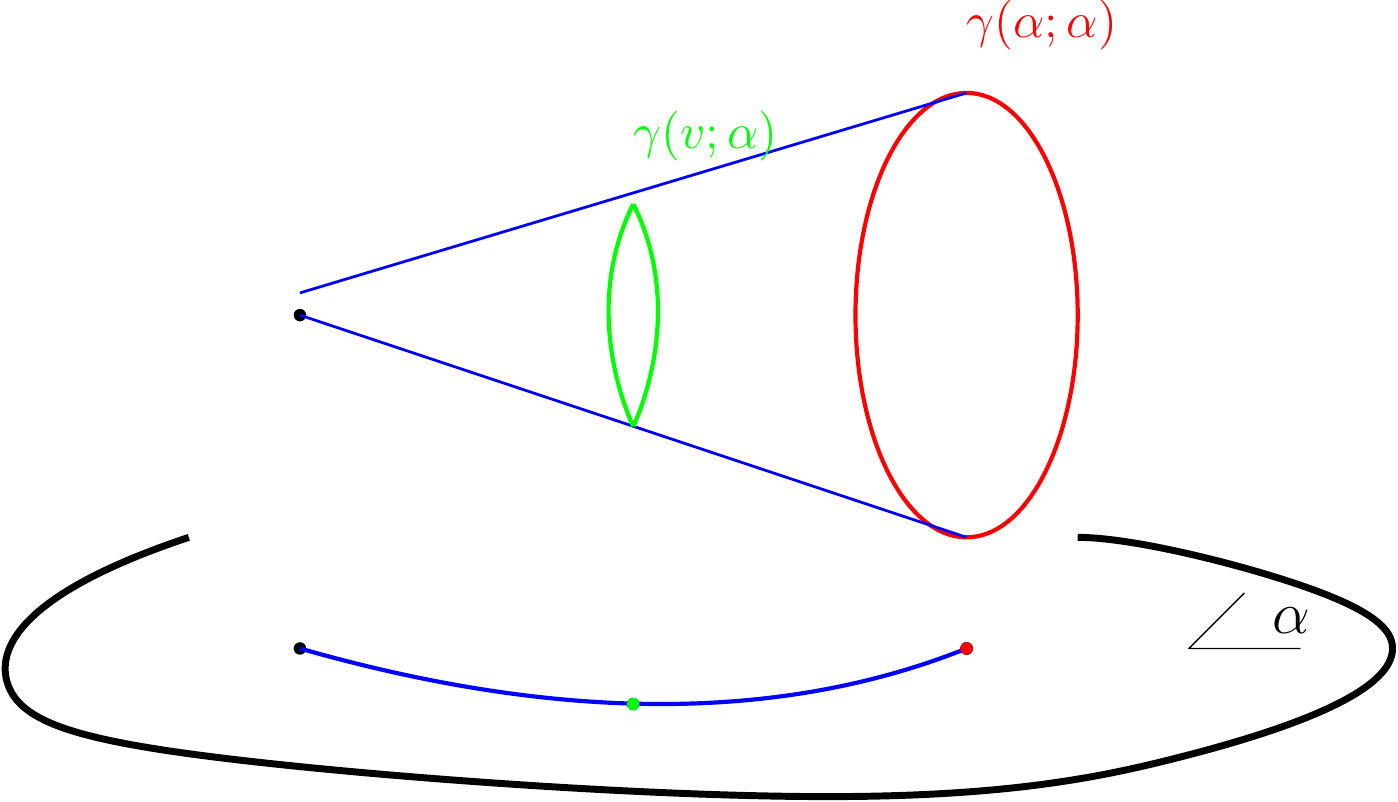}
 \caption{Region ``swept out'' by singular cycles.} \label{figurethimble}
\end{figure}

Again the orbifold point $\psi=0$ plays a special role, it is the only point $v$ in the moduli space where $\int_{\gamma(v;\alpha)}\Omega(v)=0$ for any $\alpha$, the vanishing of the integral is resulted from that of the holomorphic top form $\Omega$. Hence if we take this point as the reference point, then the quantity~$J(\alpha)$ is the area of the holomorphic form ${dv\over v}\wedge \Omega(v)$ of the cylinder swept out by these singular cycles. One can move the $\psi$ factor in~$\Omega$ to the cycle part. Then the singular cycle vanishes at the orbifold point and the cylinder becomes a disk. See Fig.~\ref{figurethimble} for an illustration.

In this way, the extra solution captures the global information of the family, as opposed to the normalized period integrals which can be def\/ined locally in the family and which does not rely on the global structure. Note that the reference point can be taken to be any point in the base of the family, the resulting chain integral carry the same amount of information through the singular cycles (and also the Beltrami dif\/ferential below), due to the algebraicity of the family.

Alternatively, the quantity
\begin{gather*}
\int_{\gamma(v;\alpha)}\Omega(v)= (u_{1}(v)u_{2}(\alpha)-u_{1}(\alpha)u_{2}(v)=\int_{\mathcal{E}_{\alpha}}\Omega(v)\wedge\Omega(\alpha)
\end{gather*}
measures the deviation of the two complex structures corresponding to $\Omega(v)$, $\Omega(\alpha)$ determined through the Torelli theorem. More precisely, one can parametrize $\Omega(v)$ in terms of the Beltrami dif\/ferential (in a suitable trivialization $\Omega(\alpha)$, $\Omega^{*}(\alpha)$ of~$H^{1}(\mathcal{E}_{\alpha},\mathbb{C})$ such that $\int_{\mathcal{E}_{\alpha}}\Omega(\alpha)\wedge \Omega^{*}(\alpha)=1$) by
\begin{gather*}
\Omega(v)=h(v;\alpha)(\Omega(\alpha)-\mu(v;\alpha)\Omega^{*}(\alpha)),
\end{gather*}
where $h(v;\alpha), \mu(v;\alpha)$ are holomorphic in $v$ but not in $\alpha$.
Then one has
\begin{gather*}
J(\alpha)=\int^{\alpha} {1\over v}dv\int_{\gamma(v;\alpha)}\Omega(v) =\int^{\alpha} {1\over v}h(v;\alpha) \mu(v;\alpha)dv .
\end{gather*}

\begin{Remark}
We can do a local calculation as follows. Fixing a choice of the section $\Omega$, we can write $\Omega(v)=\omega_{0}(v)dz_{v}$, where $dz_{v}$ is the complex coordinate on the universal cover of the elliptic curve $\mathcal{E}_{v}\cong \mathbb{C}/(\mathbb{Z}\oplus \mathbb{Z}\tau(v))$. By choosing a marking on the (generic) reference f\/iber $\mathcal{E}_{\alpha}$, the Beltrami dif\/ferential is given by the Cayley transform through
\begin{gather*}
dz_{v}=h(v;\alpha)(dz_{\alpha}-\mu(v;\alpha)(d\bar{z}_{\alpha})), \qquad
h(v;\alpha)={\tau(v)-\overline{\tau(\alpha)}\over \tau(\alpha)-\overline{\tau(\alpha)}},\qquad
\mu(v;\alpha)={\tau(v)-\tau(\alpha)\over \tau(v)-\overline{\tau(\alpha)}}.
\end{gather*}
It follows that, as already computed from the Wronskian method,
\begin{gather*}
\int_{\gamma(v;\alpha)}\Omega(v)=-(\tau(v)-\tau(\alpha))\omega_{0}(v)\omega_{0}(\alpha).
\end{gather*}
\end{Remark}

As pointed out above, the orbifold singularity point has the special property that there are two linearly independent vanishing periods corresponding to~$\pi_{1}(\alpha)$, $\pi_{2}(\alpha)$ in \eqref{eqnperiods}, while for a~generic point $\alpha$ one has only one cycle such that~\eqref{eqnvanishingperiod} is satisf\/ied.

The limit of the singular cycles at the orbifold point can be computed directly through the period calculation as follows. Since the singular cycle is independent of the marking, for computations we take~$A$,~$B$ to be the monodromy invariant cycles at the inf\/inity cusp and zero cusp respectively. There is no ambiguity in~$A$,~$B$ at the two cusps respectively, but they of course suf\/fer non-trivial monodromies elsewhere.
Their period integrals are as displayed in~\eqref{eqnanalyticcontinuation}
\begin{gather*}
\omega_{0}(\alpha)=\tilde{\pi}_{1}(\alpha)+\tilde{\pi}_{2}(\alpha),\qquad
\omega_{1}(\alpha)=-\rho \kappa\tilde{\pi}_{1}(\alpha)+\rho^2 \kappa\tilde{\pi}_{2}(\alpha).
\end{gather*}
It follows that near the orbifold point $\alpha=\infty$ or equivalently $\psi=0$ (here $\kappa=i/\sqrt{3}$)
\begin{gather*}
\gamma(\alpha;\alpha)=A\int_{B}\Omega(\alpha)-B\int_{A}\Omega(\alpha) =\tilde{\pi}_{1}(\alpha) (-\rho \kappa A-B)+\tilde{\pi}_{2}(\alpha) \big(\rho^2 \kappa A-B\big).
\end{gather*}
It is this particular linear combination of cycles that the nearby singular cycles converge to at the orbifold singularity.

\subsection[Chains on the elliptic curves and orbifold singularities on the moduli space]{Chains on the elliptic curves and orbifold singularities\\ on the moduli space}

In this section, we shall f\/ind a chain $C(\psi)$ on the elliptic curve $\mathcal{E}_{\psi}$ so that the resulting integral $\int_{C(\psi)}\Omega(\psi)$
gives the extra solution to the GKZ system other than the period integrals.

\subsubsection{Weierstrass model}\label{secWeierstrassmodel}

We f\/irst motivate the discussion by reviewing the well-studied example of Weierstrass family of elliptic curves.

As mentioned above in Section~\ref{secintro}, the derivation of dif\/ferential operators from GKZ symmetries can be applied to any family of algebraic varieties. In particular, we can apply the same discussion to the Weierstrass family
\begin{gather*}
Y^{2}=4X^3-g_{2}X-g_{3}.
\end{gather*}
The GKZ operator is computed to be
\begin{gather*}
\theta_{w}\left(\theta_{w}-{1\over 4}\right)\left(\theta_{w}-{1\over 2}\right)
-w\left(\theta_{w}+{3\over 4}\right) \left(\theta_{w}+{1\over 12}\right) \left(\theta_{w}+{5\over 12}\right), \\ w=1-{1728\over j}=27 {g_{3}^{2}\over g_{2}^3}.
\end{gather*}
The discussion by \cite{Duke:2008, Eichler:1982} implies that the extra solution is provided by the following chain integral
\begin{gather}\label{eqnWeierstrasschain}
\int_{0}^{\infty}{dX\over Y}=\int_{[\mathcal{P}(z_{0}), \mathcal{P}'(z_{0}),1]}^{[\mathcal{P}(0), \mathcal{P}'(0),1]}{dX\over Y},
\end{gather}
where $z_{0}$ is such that $\pm z_{0}$ are zeros of the Weierstrass $\mathcal{P}$-function. When pulled back to the complex $z$-plane (as the universal cover of the elliptic curve) via the Weierstrass embedding, the extra solution above is half of the chain integral{\samepage
\begin{gather}\label{eqnWeierstrasschainzplane}
\int_{-z_{0}}^{z_{0}}dz,
\end{gather}
where $dz$ is the standard holomorphic top form on the complex $z$-plane.}

This singles out the special role of the point determined by $g_{3}=0$ corresponding to the orbifold point $w=0$ in the moduli space, at which the chain $z_{0}-(-z_{0})$ on the complex $z$-plane vanishes. Intuitively, what is happening is that if one thinks of the elliptic curve as a $2:1$ cover over the $x$-plane, the chain integral mentioned above measures the ``distance'' of the two covering sheets. Its vanishing does not create a change in the topology as it does not lead to a~singular curve. However, since the value of~$g_{3}$ is non-vanishing nearby but is vanishing at the orbifold point, the vanishing of the chain integral does ref\/lect a change in the complex structure.

\begin{Remark}
This chain integral is actually the Abel--Jacobi map attached to the divisor $q-p$ given by the above two points. It appears in the study of the mixed Hodge structure of the singular curve
\begin{gather*}
Y^2=\big(4X^2-g_{2}X-g_{3}\big)X^2,
\end{gather*}
and is the obstruction to the isomorphism between the mixed Hodge structure of this singular curve and the Hodge structure of its normalization, i.e., the Weierstrass curve. Furthermore, it is the limit of a period for the genus two curve obtained by deforming the above singular curve. See~\cite{Carlson:2003} for a nice account of discussions on this.
\end{Remark}

A more natural way to look at the extra solution is to expand the corresponding oscillating integral around the orbifold point $w=\infty$ or equivalently $g_{2}=0$. The reason is that when performing the oscillating integral one is led to the following procedure (by applying change of variables)
\begin{gather*}
\int e^{-Y^2 Z+4X^3-g_{2}X Z^2-g_{3} Z^3}dXdYdZ\\
\qquad {} \implies \int Z^{-{1\over 2}} e^{-Y^2 }e^{4X^3} e^{-g_{2}X Z^2} e^{-g_{3} Z^3}dXdYdZ\\
\qquad{} \implies \int Z^{-{1\over 2}} g_{3}^{-{1\over 3}}e^{-Y^2 }e^{4X^3} e^{- Z^3}\sum_{k=0}^{\infty} {\big({-}g_{2}g_{3}^{-{2\over 3}} XZ^2\big)^{k}\over k!}dXdYdZ.
\end{gather*}
Here the integral contour needs to be chosen appropriately to deal with the convergence issue. As one can see, evaluating the integral not only provides series solutions to the GKZ system, but also picks out the natural coordinate for the expansion. From the discussion about the gauged symmetries in~\eqref{eqngaugedsymHesse}, it is easy to see that the orbifold point is always singled out according to where the polynomial~$F$ becomes a~Fermat type under a suitable coordinate change. Also the degree of the~$\mathcal{D}_{\mathrm{GKZ}}$-operator can be read of\/f easily from the action of the gauged symmetries which in particular induces action on the integral contours. This is in agreement with the result obtained by examining the linear relations in the toric data for hypersurfaces in toric varieties for example.

Therefore, both to see the gauged symmetries and to match the expansion parameter, we think of the Weierstrass elliptic curve as a $3:1$ cover over the $Y$-plane. For a generic member, there are four simply-branched points determined by, setting $f(X,Y)=Y^2-(4X^3-g_{2}X-g_{3})$,
\begin{gather*}
f=0,\qquad f_{X}=0,
\end{gather*}
as well as one 2-branched point at $Y=\infty$. The Deck group action (or the Galois action) on the covering sheets gets enhanced exactly when the simply-branched points collide. That is, when the system
\begin{gather*}
f=0,\qquad f_{X}=0,\qquad f_{XX}=0.
\end{gather*}
has non-trivial solutions. This is possible exactly at the orbifold point where $g_{2}=0$. The four simply-branched points now become two 2-branched points $Y_{o}$ determined by
\begin{gather*}
Y_{o}^2=-g_{3}.
\end{gather*}
Now we can reinterpret the chain integral in \eqref{eqnWeierstrasschain} or \eqref{eqnWeierstrasschainzplane} as the following one, which is naturally def\/ined on the $3:1$ covering,
\begin{gather*}
-\int_{-Y_{o}}^{Y_{o}}{dY\over f_{X}},
\end{gather*}
where the integral contour above means any sheet covering a path connecting the two points $-Y_{o}$,~$Y_{o}$ which may or may not pass through the branch points.

Note that carrying out the same consideration to the $2:1$ cover over the $X$-plane can only see the cusp singularities other than the orbifold singularity. As a result one can not get the full action of the gauged symmetries from the $2:1$ cover picture.

\subsubsection{Hesse pencil}

The above discussion suggests that the oscillating integral sees the f\/inest possible information of the gauged symmetries by exhibiting the most possible solutions with dif\/ferent monodromy behaviors. They are ref\/lected via the Galois symmetries of the covering with the highest possible degree.

By analogy, for the Hesse pencil, we look at the orbifold singularities in the moduli space where the conf\/iguration of branch points changes. A natural candidate for the chain whose integral gives rise to the extra solution to the GKZ system would then be the path connecting points which are not branch points for generic values of the modulus but become so at the orbifold point.

We regard a generic member of the Hesse pencil as a $3:1$ cover over the $x$-plane, in the af\/f\/ine patch $z=1$. There are $6$ simply-branched points $x_{b}$ determined by the equation
\begin{gather}\label{eqnbranchvariety}
\big(x_{b}^3+1\big)^3=4\psi^3 x_{b}^3.
\end{gather}
We denote the $6$ solutions by
\begin{gather*}
x_{b,1}, \ x_{b,2}=\rho x_{1},\qquad x_{b,3}=\rho^2 x_{1},\qquad
x_{b,4}={1\over x_{1}} , \qquad x_{b,5}={1\over x_{2}} ,\qquad x_{b,6}={1\over x_{3}} .
\end{gather*}

The symmetry of the elliptic curve for a generic value of $\psi$ given in \eqref{eqntranslationaction} is closely related to Galois symmetry of~\eqref{eqnbranchvariety}. To be more precise, the action~$\sigma_{1}$ in~\eqref{eqntranslationaction} induces
\begin{gather*}
\gamma_{3}\colon \ (x,y)\mapsto (\rho x,\rho^2 y).
\end{gather*}
The $\mathbb{Z}_{2}$ action on the complex plane as the universal cover of the elliptic curve induces $(x,y)\mapsto (y,x)$. Combing this with the $\sigma_{2}$ action in \eqref{eqntranslationaction}, one gets a symmetry of the covering
\begin{gather*}
\gamma_{2}\colon \ (x,y)\mapsto \left({1\over x},{y\over x}\right).
\end{gather*}
These actions are the Galois symmetries $\boldsymbol{\mu}_{3}\times \boldsymbol{\mu}_{2}$~\eqref{eqnbranchvariety} def\/ining the branch variety (not the Deck group transform of the $3:1$ covering).

Above a branch point $x_{b}$, the covering has three sheets determined through
\begin{gather*}
y^{3}-3\psi x_{b}y+\big(x_{b}^3+1\big)=(y-y_{b})^2 (y+2y_{b}), \qquad y_{b}^{2}=\psi x_{b}.
\end{gather*}
The Deck group of the covering is $\boldsymbol{\mu}_{2}$. This group gets enhanced at the point $\psi=0$, with $y_{b}=0$. It contains the further symmetry
\begin{gather}
(x,y,z)\mapsto (x, \rho y,z),
\end{gather}
which only exists on the f\/iber corresponding to orbifold singularity $\psi=0$ (where the elliptic curve has the extra symmetry). Note that there are other points $\psi^3=1,\infty$ such that the branch variety~\eqref{eqnbranchvariety} degenerates, but at these points the elliptic curve f\/ibers become singular and the $3:1$ coverings are only rational maps.

\begin{figure}[h]\centering
\includegraphics[scale=0.6]{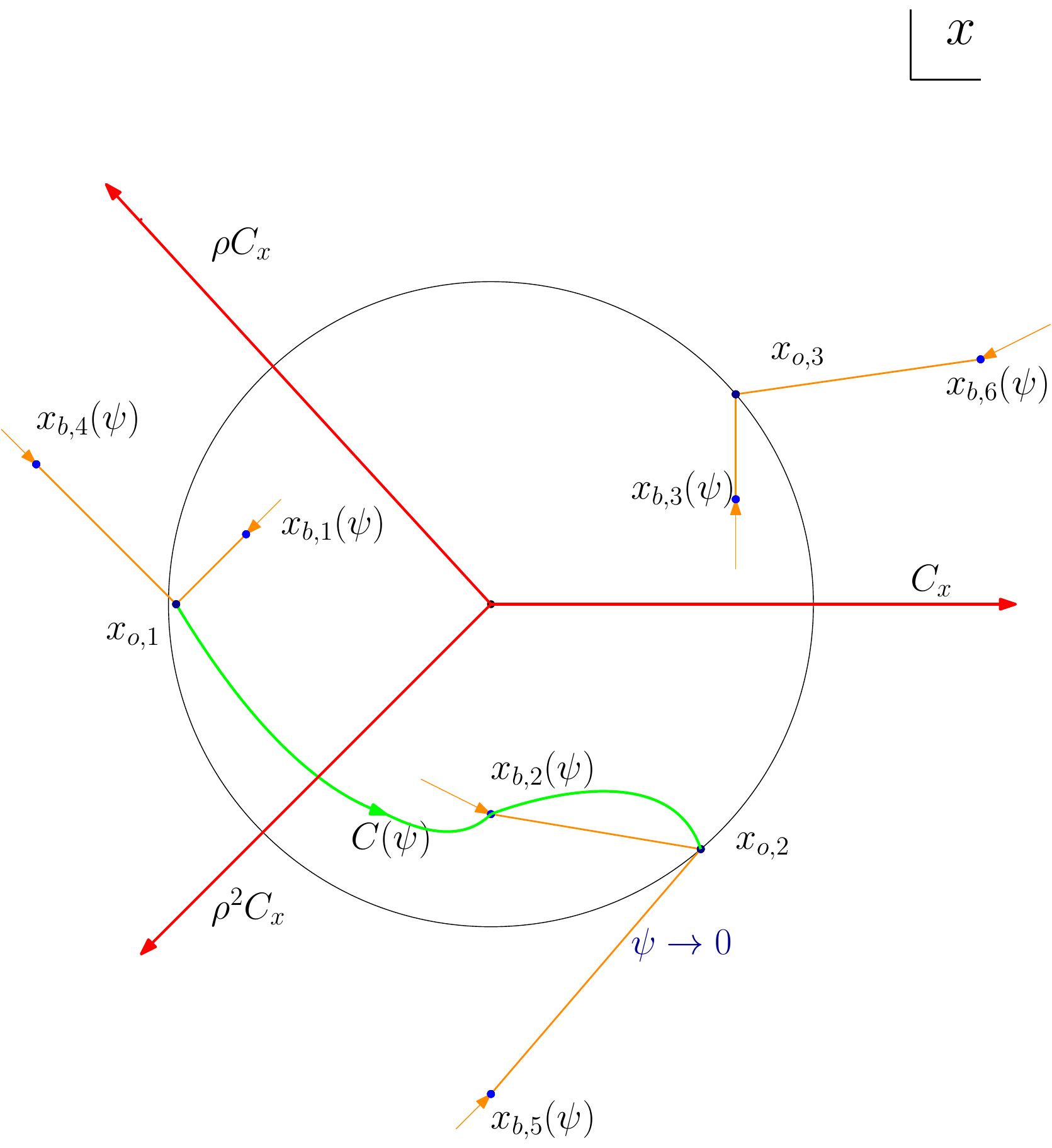}
\caption{Branch conf\/iguration of the $3:1$ cover. As $\psi\rightarrow 0$, the branch points $x_{b,k}$, $x_{b,k+3}$ collide to $x_{o,k}=-\rho^{k}$, $k=1,2,3$.} \label{figurebranch}
\end{figure}

It is easy to see that the branch points given by $x_{b,k}$, $x_{b,k+3}$ collide and gives $x_{o,k}=-\rho^{k}$, $l=1,2,3$ at the orbifold point $\psi=0$. Now on a generic f\/iber, above the point~$x_{o,k}$, the corresponding $y$-values of the points on the elliptic curve satisfy
\begin{gather*}
y_{o,k}^{3}-3\psi x_{o,k}y_{o,k}+\big(x_{o}^3+1\big)=y_{o,k}^{3}-3\psi x_{o,k}y_{o,k}=0.
\end{gather*}
The solution $y_{o,k}=0$ gives a 3-torsion point on the elliptic curve. The other two solutions satisfy $y_{o,k}^{2}=3\psi x_{o,k}$. See Fig.~\ref{figurebranch} for an illustration of the degeneration of the branch variety.

Now for a generic value of $\psi$, we take a path $C(\psi)$ on the elliptic curve with endpoints $[x_{o,k},(3\psi x_{o,k})^{1\over 2},1]$, $[x_{o,l},0,1]$, as depicted in Fig.~\ref{figurebranch}. Here $k$ could be the same as $l$.

\begin{figure}[h]\centering
\includegraphics[scale=0.66]{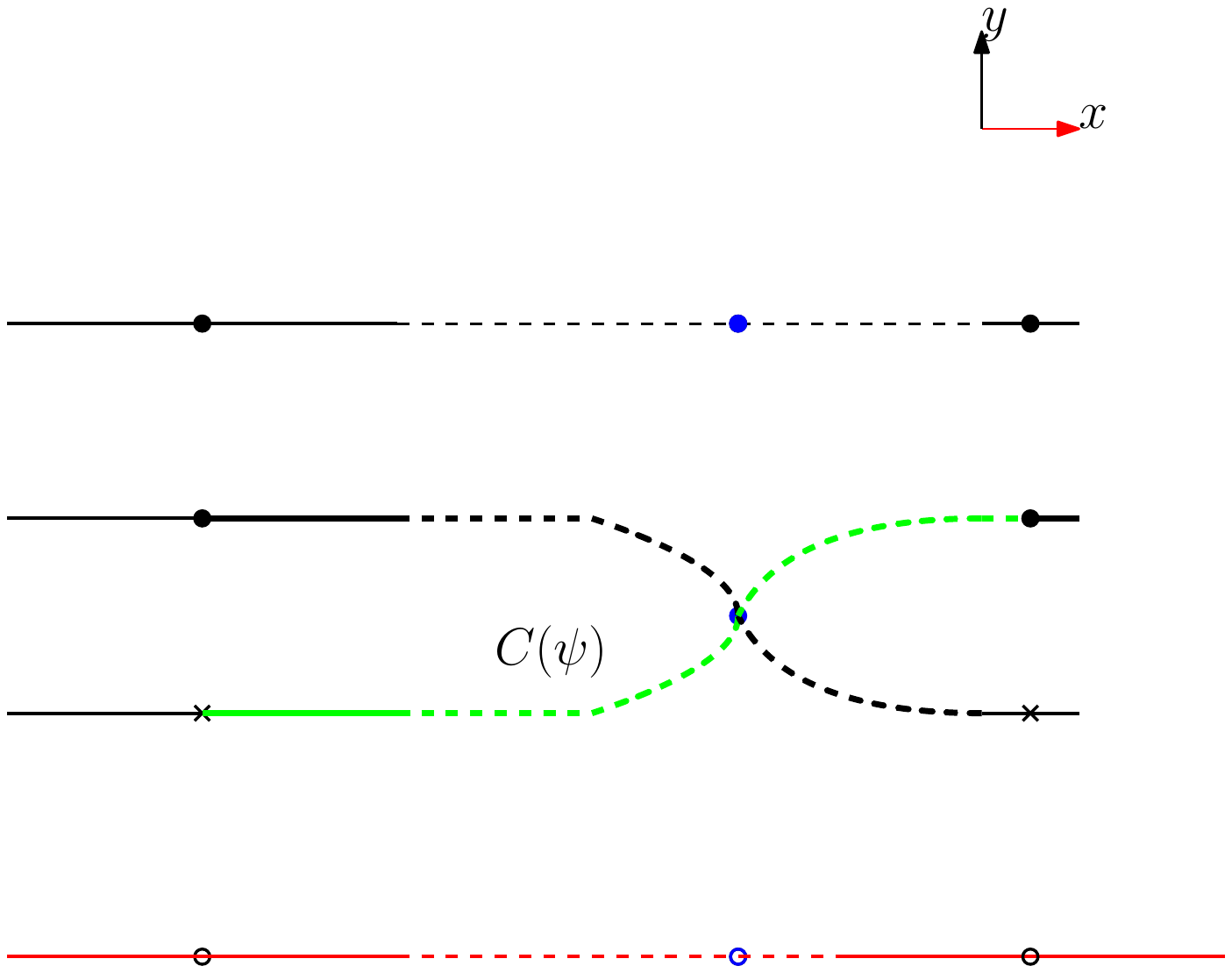}
 \caption{Chain on the elliptic curve passing through a branch point.} \label{figure1dchain}
\end{figure}

Note that a path like this must pass through a branch point. The dif\/ference between two such paths are cycles and hence their integrals are dif\/fered by period integrals. See Fig.~\ref{figure1dchain} for an illustration.

Note also that if instead one takes a path with endpoints $[x_{o,k},0,1]$, $[x_{o,l},0,1]$, $k\neq l$, then one gets a~chain connecting two 3-torsion points and the resulting integral is a period integral as mentioned before in Section~\ref{secdifferencesofchainintegrals}.

We consider the chain integral
\begin{gather}\label{eqn1dchainintegrald}
K(\psi):=\int_{C(\psi)}{\psi dx\over 3y^2-3\psi x}.
\end{gather}
We then have the following result.
\begin{Theorem}\label{thmchainintegral}
The integral \eqref{eqn1dchainintegrald} gives a solution to the GKZ system for the Hesse pencil and is not a period integral.
\end{Theorem}
\begin{proof}
From the Grif\/f\/iths--Dwork method, we can f\/ind the exact term in the Picard--Fuchs operator acting on the holomorphic top form to be
\begin{gather*}
 \big( \theta_{\psi}^2- \psi^{-3} (\theta_{\psi}-2)(\theta_{\psi}-1) \big)\left( {\psi dx\over 3y^2-3\psi x} \right)
=\psi^{-3}\circ \mathcal{L}_{\mathrm{PF}}\left( {\psi dx\over 3y^2-3\psi x} \right) =d\left({\psi x\over f_{y}}\right).
\end{gather*}
Note that by construction the endpoints of the path $C(\psi)$, when parametrized by~$x$, are locally constant and hence annihilated by the derivatives.
Then by Stokes theorem, one can immediately check the in-homogeneous Picard--Fuchs equation
\begin{gather*}
\big(\psi^{-3}\circ \mathcal{L}_{\mathrm{PF}}\big) K(\psi) =\int_{C(\psi)} d\left({\psi x\over y}\right)=\left({\psi x\over f_{y}}\right)\Big|_{\partial C(\psi) }={1\over 6}-\left(-{1\over 9}\right) \neq 0.
\end{gather*}
Therefore, these chain integrals do give solutions to the GKZ operator $\mathcal{D}_{\mathrm{GKZ}}=\theta\circ (\psi^{-3}\circ \mathcal{L}_{\mathrm{PF}})$ which are not periods.
\end{proof}

\begin{Remark} The endpoints of the chain belong to $\mathcal{E}_{\psi}\cap \{x^3+z^3=0\}$. Under the uniformization by he $\theta$-functions, see, e.g.,~\cite{Dolgachev:1997}, these end points are zeros of certain theta functions and carry interesting arithmetic information, as is the case in the Weierstrass family discussed in Section~\ref{secWeierstrassmodel}. In terms of the toric coordinates corresponding to the characters, this is similar to the situation that appears in open mirror symmetry \cite{Lerche:2003hs, Lerche:2001cw, Lerche:2002ck, Lerche:2002yw, Mayr:2001} which again seems to indicate that the chain integral is related to the enumerative geometry in the A-model under mirror symmetry.
\end{Remark}

\begin{Remark}
One can also consider the higher Frobenius functions appearing as the coef\/f\/i\-cients in the $\epsilon$-expansion of the function
\begin{gather*}
\omega_{0}(\alpha,\epsilon)=\sum_{n=0}^{\infty} Q(\epsilon){\Gamma(3n+3\epsilon+1)\over \Gamma(n+\epsilon+1)^3} \left({\alpha\over 3^{3}}\right)^{n+\epsilon}:=\sum_{k=0}^{\infty} f_{k}(\alpha)\epsilon^{k},\qquad Q(\epsilon)={ \Gamma(\epsilon+1)^{3} \over \Gamma(3\epsilon+1)}.
\end{gather*}
This is the deformation of the period in~\eqref{eqnanalyticcontinuation} when applying the Frobenius method to solve for the solutions to the Picard--Fuchs equation. It satisf\/ies, recall~\eqref{eqnPFinalphacoordinate},
\begin{gather*}
 \tilde{\mathcal{L}}_{\mathrm{PF}} \omega_{0}(\alpha,\epsilon)=\left({\alpha\over 3^{3}}\right)^{\epsilon}\epsilon^{2},\qquad
 \tilde{\mathcal{D}}_{\mathrm{GKZ}}\omega_{0}(\alpha,\epsilon)=\left({\alpha\over 3^{3}}\right)^{\epsilon}\epsilon^{3}.
\end{gather*}
Therefore, one has
\begin{gather*}
 \tilde{\mathcal{L}}_{\mathrm{PF}} f_{k}(\alpha)= {(\ln {\alpha\over 3^{3}} )^{k-2}\over (k-2)!}, \qquad
 \tilde{\mathcal{D}}_{\mathrm{GKZ}} f_{k}(\alpha)={(\ln {\alpha\over 3^{3}} )^{k-3}\over (k-3)!}.
\end{gather*}
Here we have used the convention that negative powers of $(\ln {\alpha\over 3^{3}} )$ give zero. Besides $\{f_{0},f_{1}\}$ which are period integrals and
$\{f_{0},f_{1},f_{2}\}$ which are chain integrals, the higher Frobenius functions $\{f_{k},\, k\geq 3\}$, which can be solved by the Wronskian method as in Section~\ref{secWronskian}, are also interesting on their own. For example, they carry interesting arithmetic meanings, corresponding to the counting of rational points of the Hesse elliptic curves \cite{Candelas:2000, Candelas:2003}. Furthermore, when one regards the variable~$\epsilon$ as the hyperplane class of $\mathbb{P}^{2}$, then $\omega_{0}(\alpha,\epsilon)$ gives Givental's (twisted) $I$-function valued in the cohomology ring and the factor~$Q(\epsilon)$ is the $\Gamma$-class which appears naturally in mirror symmetry, see~\cite{Galkin:2016, Golyshev:2016, Hosono:1995, Libgober:1999}. After passing to the equivariant cohomology corresponding to the diagonal torus action, the higher Frobenius functions then appear as the coef\/f\/icients of the equivariant version of the $I$-function expanded in the equivariant parameter. We wish to discuss their geometric meanings in a future work.
\end{Remark}

\subsubsection{Legendre family}

We conclude this section with some discussions on the Legendre family whose af\/f\/ine equation is
\begin{gather}\label{eqnLegendre}
y^{2}=x(x-1)(x-\lambda),\qquad j(\lambda)=2^8 {(\lambda^{2}-\lambda+1)^{3}\over \lambda^2(\lambda-1)^2}.
\end{gather}
From the derivation of GKZ system using the GKZ symmetries, we can see that the operator $\mathcal{D}_{\mathrm{GKZ}}$ is a 2nd order operator and hence coincides with the Picard--Fuchs operator. This also agrees with the earlier discussion on the relation between the extra solution and orbifold singularities. Namely, in this case one can check that the Picard--Fuchs equation has no orbifold singularity. One can also see this by using the standard fact that the base of the family is parametrized by the modular curve $\Gamma(2)\backslash\mathcal{H}^{*}$ which has no elliptic f\/ixed point.

However, from the evaluation of the oscillating integral, one can see that
\begin{gather*}
\begin{split}
& \int e^{-(y^2 z-(x^3-(\lambda+1)x^2 z+\lambda xz^2))}dxdydz\\
& \qquad{} \sim
\lambda^{-{1\over 2}} \int e^{-y^2} e^{x^3}e^{z^2}\sum_{k=0}^{\infty} \left({\lambda+1\over \lambda^{1\over 2}}\right)^{k}
{x^{{3\over 2}k-{1\over 2}}z^{k-{1\over 2}} \over k!}dxdydz.
\end{split}
\end{gather*}
Therefore, the oscillating integral naturally singles out the coordinate $ \alpha=({\lambda+1)/ \lambda^{1\over 2}}$ for the expansion parameter.\footnote{The transformation from the $\lambda$-parameter to this $\alpha$-parameter is induced by a $2$-isogeny, as can be seen through the elliptic $\kappa$-modulus.} Moreover, the gauged symmetries would give rise to at least $6$ solutions with dif\/ferent monodromy behaviors. This is however not a contradiction to the statement that the $\mathcal{D}_{\mathrm{GKZ}}$ is of second order. The reason is that in deriving the $\mathcal{D}_{\mathrm{GKZ}}$ using the GKZ symmetries, only scalings on the parameter $\lambda$ are allowed and hence those act by scalings on the new parameter $\alpha$ is not included.

Furthermore, the locus at which $\alpha=0$ corresponds to the point $\lambda=-1$ or $j=1728$ according to the formula for the $j$-invariant in~\eqref{eqnLegendre}. Hence indeed when parametrized by $\alpha$ the above expansion of the oscillating integral occurs near an orbifold point.

One can again obtain chain integrals by studying the branch conf\/iguration of the $3:1$ cover realization for the elliptic curve. Now the enhancement of the Galois symmetry takes place at $\lambda=-\rho,-\rho^2$ where $j=0$. This corresponds to a dif\/ferent way of performing the oscillating the integral above by f\/irst applying the following change of variables and then evaluating
\begin{gather*}
-y^2 z+\big(x^3-(\lambda+1)xz^2+\lambda z^3\big)=-y^2 z+\left(x-{\lambda+1\over 3}\right)^3+\left({\lambda+1\over 3}\right)^{3}z^3\\
\qquad{}
-\left( {\lambda^2-\lambda+1\over 3}\right)\left(x-{\lambda+1\over 3}\right)z^2- {(\lambda+1)(\lambda^2-\lambda+1)\over 3^2 } z^2.
\end{gather*}

In summary, the oscillating integrals of\/fer more than the solutions to the GKZ systems. It has the f\/inest information about the gauged symmetries which includes the GKZ scaling symmetries as a subset.

\section{Period integrals in a compact Calabi--Yau threefold}\label{secembedding}

In this section, we shall explain the relation between the two dif\/ferential operators, namely $\mathcal{D}_{\mathrm{GKZ}}$ in \eqref{eqnrightfactorization} and $\mathcal{L}_{\mathrm{CY}_{3}}$ in \eqref{eqnleftfactorization}. We shall see that they are dif\/ferent pieces of the same Picard--Fuchs system of a compact Calabi--Yau threefold.

Recall that in Section \ref{seclocalCYintegrals} we explained that the $3$d oscillating integrals and $2$d real integrals in \eqref{eqndimreduction} should be interpreted as ones on a non-compact Calabi--Yau. We now push this idea further.

We f\/irst note that the members in the Hesse pencil correspond to the sections of the anti-canonical divisor of the toric variety~$P$ whose polytope is generated by $(1,0)$, $(0,1)$, $(-1,-1)$. This polytope is a ref\/lective polytope and def\/ines the following toric variety
\begin{gather}\label{eqnmirrorP2}
P=\mathbb{P}^{2}/G, \qquad G=\big\{(\rho^{n_{1}},\rho^{n_{2}},\rho^{n_{3}})\,|\,n_{1}+n_{2}+n_{3}=0~\mathrm{mod}~3\big\}.
\end{gather}
The invariants of $G$ are the monomials $x_{1}^3$, $x_{2}^3$, $x_{3}^{3}$, $x_{1}x_{2}x_{3}$ among all cubic monomials. The induced action on a generic member of the Hesse pencil is the one generated by~$\sigma_{1}$ in~\eqref{eqntranslationaction}. The quotient therefore gives the $3$-isogeny of the Hesse pencil and is the mirror of the Hesse pencil according to~\cite{Batyrev:1994hm}. It can be checked by using the GKZ symmetries or by evaluating the oscillating integrals that these two elliptic curve families share the same GKZ operators. Hence for the purpose of studying the solutions to the GKZ system, there is no dif\/ference between these two families. See~\cite{Zhou:2016} for discussions on these facts and some arithmetic aspects of the mirror symmetry.

Now the quotient of the Hesse pencil is naturally interpreted as sections of the canonical bundle $K_{P}$ of~$P$. There is a natural compactif\/ication \cite{Candelas:1993dm, Chiang:1999tz} $X$ of $K_{P}$. A certain limit of $X$ gives rise to the variety $K_{P}$, as the mirror of $K_{\mathbb{P}^{2}}$, whose Picard--Fuchs operator is displayed in~\eqref{eqnleftfactorization}. This idea is used frequently in the literature to study mirror symmetry for non-compact Calabi--Yau manifolds. As we shall review in Section \ref{secreviewcompactification}, this compactif\/ication also encodes the full information of the quotient of the Hesse pencil by $G$, as the mirror of the Hesse pencil, including the GKZ operator $\mathcal{D}_{\mathrm{GKZ}}$ in~\eqref{eqnrightfactorization}.

It is then natural to expect a relation between the two geometries--mirror of Hesse pencil and mirror of $K_{\mathbb{P}^{2}}$~-- by embedding them in the same ambient space~$X$. The properties about the GKZ/Picard--Fuchs system should be independent of the choice for the compactif\/ication though.

\subsection{Review of the compactif\/ication}\label{secreviewcompactification}

We now recall the construction of the compactif\/ication following \cite{Candelas:1993dm, Chiang:1999tz}. The A-model is an elliptic f\/ibration over $\mathbb{P}^{2}$. The total space $\check{X}$ is a Calabi--Yau hypersurface in a~toric variety. For the mirror geometry $X$, the toric data gives the family of varieties $X$ whose Zariski open sets are described by the equation
\begin{gather*}
\Xi=b_{0}+Z_{3}^2 Z_{4}^{3} \big(a_{1}Z_{1}+a_{2}Z_{2}+a_{3} Z_{1}^{-1}Z_{2}^{-1}+a_{0}\big) +a_{4}Z_{3}^{-1}+a_{5}Z_{4}^{-1}.
\end{gather*}
Switching to the homogeneous coordinates, this is
\begin{gather}
\Xi=\left(\prod x_{i}^{-1} \right)\big(b_{0}x_{1}x_{2}x_{3}x_{4}x_{5}+ a_{1}x_{1}^{18}+a_{2}x_{2}^{18}+a_{3}x_{3}^{18}+a_{0}x_{1}^6 x_{2}^6 x_{3}^6 + a_{4}x_{4}^3+a_{5}x_{5}^2\big)\nonumber\\
\hphantom{\Xi}{} :=\left(\prod x_{i}^{-1} \right)\xi.\label{eqnmirrorfamily}
\end{gather}
It is an elliptic f\/ibration over the base $P$ which is parametrized by $x_{1}$, $x_{2}$, $x_{3}$. We ignore the subtitles about the group actions involved which do not af\/fect the Picard--Fuchs systems we are interested in. Thinking of~$X$ as a Weierstrass f\/ibration over~$P$, we then get the identif\/ication for the divisor classes
\begin{gather}
\big(x_{i}^3=0\big), (x_{1}x_{2}x_{3}=0)=\mathcal{O}_{WP}(1),\qquad i=1,2,3,\nonumber\\
(x_{4}=0)=\mathcal{O}_{WP}(2)\otimes K_{P}^{-2},\qquad (x_{5}=0)=\mathcal{O}_{WP}(3)\otimes K_{P}^{-3},\label{eqnBmodelWeierstrassform}
\end{gather}
where $WP$ denotes the weighted projective space $\mathrm{WP}[1,2,3]$ in which the elliptic curve f\/ibers sit. This implies that the coef\/f\/icients transform as sections of certain tensor powers of $K_{P}$: the variables $a_{4}$, $a_{5}$, $b_{0}$ transform as sections of~$ K_{P}^{-1}$, while $a_{1}$, $a_{2}$, $a_{3}$, $a_{0}$ sections of~$K_{P}^{-6}$. Strictly speaking, they are sections of the corresponding relative line bundles over the base of the f\/ibration.

Note that setting $a_{4}=a_{5}=b_{0}=0$ in $\xi$ gives the equation for the Hesse pencil. The limit $b_{0}=0$ in $\xi$ gives~\cite{Chiang:1999tz} the mirror of~$K_{\mathbb{P}^{2}}$. We shall say more about this below.

By using the GKZ symmetries, one can simplify $\xi$ into
\begin{gather*}
\big(bx_{1}x_{2}x_{3}x_{4}x_{5}+ x_{1}^{18}+x_{2}^{18}+x_{3}^{18}+a x_{1}^6 x_{2}^6 x_{3}^6 +x_{4}^3+x_{5}^2\big).
\end{gather*}
Here\footnote{We have used dif\/ferent notations for the parameters from those in \cite{Candelas:1993dm}.}
\begin{gather}\label{eqncomplexmoduli}
a=(a_{1}a_{2}a_{3})^{-{1\over 3}}a_{0}, \qquad b=b_{0}(a_{1}a_{2}a_{3})^{-{1\over 18}}a_{4}^{-{1\over 3}}a_{5}^{-{1\over 2}}.
\end{gather}
There are interesting loci in the base of the family $\Xi$ parametrized by the coordinates~$(a,b)$. In particular, the point $a=b=\infty$ corresponds to the large complex structure limit. See \cite{Candelas:1993dm, Chiang:1999tz} and also \cite{Alim:2012ss, Klemm:2012sx} for details.

\subsection[Picard--Fuchs system and fundamental period of the compactif\/ied geometry]{Picard--Fuchs system and fundamental period\\ of the compactif\/ied geometry}

The period integrals are the integrals of the following form over the tubular neighborhood of cycles in $X$
\begin{gather}\label{eqnCY3integral}
\int {b_{0}\over \Xi} {dZ_{1}dZ_{2}dZ_{3}dZ_{4}\over Z_{1}Z_{2}Z_{3}Z_{4}} =\int {b_{0}\mu_{0}\over \xi}.
\end{gather}
where $\mu_{0}$ denotes the standard meromorphic $4$-form in the ambient space which has a pole of order one at inf\/inity. The Picard--Fuchs system can be derive from the GKZ symmetries and are given as follows\footnote{These Picard--Fuchs operators are derived by factoring out some dif\/ferential operators
from the left in the GKZ $\mathcal{Z}$-operators. One can also study the extra solutions to the GKZ system of the current Calabi--Yau threefold by embedding it into a variety of higher dimension, similar to what will be discussed below. But we shall not discuss them in this work.}
\begin{gather*}
\mathcal{D}_{1}= {1\over (-3)^{2}(-2)^{3}}\theta_{a_{0}} \theta_{b_{0}}-a b^{-6}(\theta_{b_{0}}-1)(\theta_{b_{0}}-5),\nonumber\\
\mathcal{D}_{2}={1\over (-18)^{3}}(\theta_{b_{0}}+6\theta_{a_{0}})^3-a^{-3}(\theta_{a_{0}}-1)(\theta_{a_{0}}-2)\theta_{a_{0}} .
\end{gather*}
In terms of the coordinates $a$, $b$, we get
\begin{gather}
\mathcal{D}_{1}= {1\over (-3)^{2}(-2)^{3}}\theta_{a} \theta_{b}-a b^{-6}(\theta_{b}-1)(\theta_{b}-5),\nonumber\\
\mathcal{D}_{2}={1\over (-18)^{3}}(\theta_{b}+6\theta_{a})^3-a^{-3}(\theta_{a}-1)(\theta_{a}-2)\theta_{a} .\label{eqnPFsystemab}
\end{gather}

The fundamental period\footnote{The unique (up to scaling) regular period near the large complex structure limit given by $a=b=\infty$.} can be obtained directly by manipulating the series expansion in~\eqref{eqnCY3integral} with a suitable choice for the integral contour, as done in~\cite{Candelas:1993dm, Chiang:1999tz}. It is given by
\begin{gather}
\omega_{0}(a,b)=\sum_{n,m=0}^{\infty}{\Gamma(18n+6m+1)\over \Gamma(9n+3m+1)\Gamma(6n+2m+1)\Gamma(n+1)^3 \Gamma(m+1)}a^{m}(b^{-6})^{3n+m}\nonumber\\
\hphantom{\omega_{0}(a,b)}{} =
\sum_{k=0}^{\infty} {\Gamma(6k+1)\over \Gamma(3k+1)\Gamma(2k+1)\Gamma(k+1)}b^{-6k}U_{k}(a):=
\sum_{k=0}^{\infty}c_{k}b^{-6k}U_{k}(a),\label{eqnperiodexpansion}
\end{gather}
where
\begin{gather*}
U_{k}(a)=a^{k}\sum_{l=0}^{[{k\over 3}]} {\Gamma(k+1)\over \Gamma(l+1)^{3}\Gamma(k-3l+1)}a^{-3l}.
\end{gather*}

The above expansion \eqref{eqnperiodexpansion} amounts to solving the Picard--Fuchs system \eqref{eqnPFsystemab} in the following way. The degree $k$-piece $c_{k}b^{-6k}U_{k}(a)$ in the sum satisf\/ies $\theta_{b}=-6k$. Hence the second equation in~\eqref{eqnPFsystemab} gives the equation for $U_{k}(a)$
\begin{gather*}
\left({1\over (-18)^{3}}(-6k+6\theta_{a})^3-a^{-3} (\theta_{a}-1)(\theta_{a}-2)\theta_{a}\right) U_{k}(a)=0.
\end{gather*}
This can be simplif\/ied into
\begin{gather*}
\left( (\theta_{a}-1)(\theta_{a}-2)\theta_{a}-a^3{6^3\over (-18)^{3}}(\theta_{a}-k)^3\right) U_{k}(a)=0.
\end{gather*}
The f\/irst equation in \eqref{eqnPFsystemab} then gives recursive relations among $\{c_{k}b^{-6k}U_{k}(a)\}_{k}$ through
\begin{gather*}
\sum_{k}(-6k)b^{-6k}c_{k}\theta_{a}U_{k}=\sum_{k}(-3)^{2} (-2)^{3} ab^{-6k-6}(6k+1)(6k+5) c_{k}U_{k}.
\end{gather*}
This is simplif\/ied into
\begin{gather*}
\theta_{a}U_{k+1}= (k+1)a U_{k}.
\end{gather*}

\subsection[Embedding of the GKZ system for the Hesse pencil and the Picard--Fuchs system for the mirror geometry of $K_{\mathbb{P}^{2}}$]{Embedding of the GKZ system for the Hesse pencil \\ and the Picard--Fuchs system for the mirror geometry of $\boldsymbol{K_{\mathbb{P}^{2}}}$}

For the purpose of getting the other solutions via the Frobenius method and doing analytic continuation, one needs to extend~\cite{Candelas:1993dm} the def\/inition of $U_{k}$ to $U_{\nu}$ for complex values of~$\nu$
\begin{gather*}
U_{\nu}(a)=a^{\nu}\sum_{l=0}^{\infty} {\Gamma(\nu+1)\over \Gamma(l+1)^{3}\Gamma(\nu-3l+1)}a^{-3l}
=a^{\nu}\,{}_{3}F_{2}\left({-\nu\over 3}, {1-\nu\over 3}, {2-\nu\over 3};1,1;a^{-3}\right).
\end{gather*}
It can be analytically continued to the orbifold $a=0$ via the Barnes integral formula \cite{Candelas:1993dm}
\begin{gather}\label{eqnanalyticcontinuationoffundamentalperiod}
U_{\nu}(a)= {3^{-1-\nu} \rho^{{\nu}\over 2}\over \Gamma(-\nu)}\sum_{n=0}^{\infty}{\Gamma({n-\nu\over 3}) \over \Gamma^{2}(1-{{n-\nu}\over 3})}{(-3\rho a)^{n}\over n!}.
\end{gather}

The recursive relation is given by
\begin{gather*}
\theta_{a}U_{\nu+1}= (\nu+1)a U_{\nu}.
\end{gather*}
 It is annihilated by the operator
\begin{gather*}
\mathcal{L}_{\nu}=\left( (\theta_{a}-1)(\theta_{a}-2)\theta_{a}-a^3{6^3\over (-18)^{3}}(\theta_{a}-\nu)^3\right).
\end{gather*}
Setting $a=-3\psi$ (which makes contact with the Hesse pencil), one gets
 \begin{gather*}
\mathcal{L}_{\nu}=\big( (\theta_{\psi}-1)(\theta_{\psi}-2)\theta_{\psi}-\psi^{3}(\theta_{\psi}-\nu)^3\big).
\end{gather*}

When $\nu=0$, this is the Picard--Fuchs operator $\mathcal{L}_{\mathrm{CY}_{3}}$ in \eqref{eqnleftfactorization}.

When $\nu=-1$, this is equivalent to the operator $\mathcal{D}_{\mathrm{GKZ}}$ in \eqref{eqnrightfactorization} and it annihilates the form~${\mu_{0}\over F}$ in~\eqref{eqnmero2form}. The solution given in~\eqref{eqnanalyticcontinuationoffundamentalperiod} is exactly the one in~\eqref{eqnoscillating} up to a~constant multiple.

In general, $\mathcal{L}_{\nu}$ annihilates
\begin{gather}\label{eqnformonambientspace}
a_{0}^{\nu+1}{\mu_{0}\over F}.
\end{gather}
As explained in Section~\ref{seclocalCYintegrals}, one should think of the parameter~$a_{0}$ (previously denoted by~$sa_{0}$ in the trivialization $\mu_{0}$) as the coordinate of the f\/iber of $K_{P}$, and hence $\nu+1$ as the degree of the form in~\eqref{eqnformonambientspace} along the f\/iber direction.

Consider the analytic continuation of the expansion \eqref{eqnperiodexpansion} to the orbifold point $b=0$, then one has~\cite{Candelas:1993dm}
\begin{gather*}
\omega_{0}(a,b)={1\over 2\pi }
\sum_{k=0}^{\infty} 2^{-{1\over 2}-2k}3^{-{1\over 2}-3k}6^{{1\over 2}+6k}
{\Gamma(k+{1\over 6}) \Gamma(k+{5\over 6}) \over \Gamma(k+1)^{2}}b^{-6k}U_{k}(a)\nonumber\\
\hphantom{\omega_{0}(a,b)}{} =\sum_{k=0}^{\infty}{\Gamma(k+{1\over 6}) \Gamma(k+{5\over 6}) \over \Gamma(k+1)^{2}} (432)^{k}b^{-6k}U_{k}(a):=
\sum_{n=0}^{\infty} d_{n}b^{n}U_{-{n\over 6}}(a).
\end{gather*}
Here $\{d_{n}\}_{n}$ are some Gamma-values whose precise values are not important in the discussion here. Then we can see that both $U_{0},U_{-1}$ appear in the fundamental period as pieces in $b$-expansion of dif\/ferent degrees. They appear naturally in the expansion around the orbifold point $b=0$
as opposed to the expansion near the point $b=\infty$ in~\eqref{eqnperiodexpansion}.

More geometrically, one can expand the dif\/ferential form in \eqref{eqnCY3integral} as follows
\begin{gather*}
{b_{0}\mu_{0}\over \xi}=b_{0}\mu_{0} \sum_{k=0}^{\infty}
{1\over \big(a_{1}x_{1}^{18}+a_{2}x_{2}^{18}+a_{3}x_{3}^{18}+a_{0} x_{1}^6 x_{2}^{6} x_{3}^6+a_{4} x_{4}^{3}+a_{5} x_{5}^{2}\big)^{k+1}}\\
\hphantom{{b_{0}\mu_{0}\over \xi}}{}\times
(-b_{0})^{k} (x_{1}x_{2}x_{3}x_{4}x_{5})^{k}.
\end{gather*}
The degree zero term in $b_{0}$ in the summation gives the holomorphic volume form in~\eqref{eqnnoncompactvolumeform}. Here we treat the prefactor $b_{0}$ of $\mu_{0}$ as an overall normalization. Taking $b_{0}=0$ is equivalent to the limit when the compact Calabi--Yau threefold~$X$ degenerates to the mirror $K_{P}$ of~$K_{\mathbb{P}^{2}}$. It is actually more convenient to see the degenerating limit in the toric coordinates $Z_{i}$ on the torus in the toric variety. One writes ${b_{0}\mu_{0}/ \xi}$ as the following, from which one recognizes~\eqref{eqnnoncompactvolumeform} easily,
\begin{gather}\label{eqntoruscooordinatesdegeneration}
{b_{0}(-b_{0})^{k}\over \big(Z_{3}^2 Z_{4}^{3} (a_{1}Z_{1}+a_{2}Z_{2}+a_{3} Z_{1}^{-1}Z_{2}^{-1}+a_{0})
+a_{4}Z_{3}^{-1}+a_{5}Z_{4}^{-1}\big)^{k+1}}{dZ_{1}dZ_{2}dZ_{3}dZ_{4}\over Z_{1}Z_{2}Z_{3}Z_{4}}.
\end{gather}

Setting $b_{0}$ to zero means that one is ef\/fectively looking at the vanishing of a section of $\mathcal{O}_{\mathrm{WP}[1,2,3]}(1)$ in the f\/iber weight projective space, that is, the divisor $(x_{1}x_{2}x_{3})=0$ according to~\eqref{eqnBmodelWeierstrassform}. This def\/ines the (unique) section of the Weierstrass f\/ibration.

Alternatively, one can write
\begin{gather}
{b_{0}\mu_{0}\over \xi} =b_{0}\mu_{0}
\sum_{k=0}^{\infty}{(-1)^{k}\over (a_{1} x_{1}^{18}+a_{2} x_{2}^{18}+a_{3} x_{3}^{18}+a_{0} x_{1}^6 x_{2}^{6} x_{3}^6)^{k+1}}\nonumber\\
\hphantom{{b_{0}\mu_{0}\over \xi} =}{}\times
\big(b_{0} x_{1}x_{2}x_{3}x_{4}x_{5}+a_{4} x_{4}^{3}+a_{5} x_{5}^{2}\big)^{k}.\label{eqnGKZlimit}
\end{gather}

Now in the limit $a_{4}=a_{5}=b_{0}=0$, one recovers the meromorphic $2$-form $\mu_{0}/F$ in \eqref{eqnmero2form} that appears in the GKZ system for the Hesse pencil. Again here we have regarded the prefactor~$b_{0}$ of~$\mu_{0}$ as an overall normalization. In fact, since in the integration, the contour that gives rise to the fundamental period is parametrized in such a way that the coordinates~$x_{4}$, $x_{5}$ take values in~$S^{1}$, the above limit on the period can be induced by the limit $a_{4}=a_{5}=0$.

In terms of the geometry, this amounts to setting $x_{4}=x_{5}=0$, which cuts out the Hesse pencil in~$X$. One can also see this by examining~\eqref{eqntoruscooordinatesdegeneration} in the coordinates $Z_{i}$, $i=1,2,3,4$ on the torus. Intuitively, the Calabi--Yau threefold $X$ admits a rational map to~$\mathrm{WP}[1,2,3]$ as an elliptic f\/ibration, the f\/ibers are the Hesse elliptic curves. The equations $x_{4}=x_{5}=0$ def\/ines the f\/iber at the singular point with stabilizer $\boldsymbol{\mu}_{6}$ in~$\mathrm{WP}[1,2,3]$.

In either case, the degree in $b_{0}$ indicates the degree $\nu$ along the f\/ibration direction of~$K_{P}$.

\subsection{Interpretation in the A-model}

While it is straightforward to see the above degeneration limits by examining the def\/ining equation of the Calabi--Yau variety~$X$, it is perhaps also helpful to study these limits in the A-model geometry.

The family in the A-model is parametrized by the space of K\"ahler structures of the variety~$\check{X}$. The latter is the resolution of singularities of a degree $18$ hypersurface $\check{X}_{0}$ in the weight projective space $\mathrm{WP}[1,1,1,6,9]$ parametrized by $x_{1}$, $x_{2}$, $x_{3}$, $x_{4}$, $x_{5}$. The singularity occurs at $x_{1}=x_{2}=x_{3}=0$. The details are worked out in~\cite{Candelas:1993dm}. We now give a brief review on the intersection theory of the geometry. One denotes the strict transform of~$(x_{1}=0)$ by~$L$, and the total transform of the divisor $(x_{1}x_{2}x_{3}=0)$ by $H=3L+E$, where~$E$ is the class of the exceptional divisor. The intersections are
\begin{gather*}
H^3=9, \qquad H^{2}L=3,\qquad HL^{2}=1,\qquad L^{3}=0.
\end{gather*}

\begin{figure}[h]\centering
\includegraphics[scale=0.6]{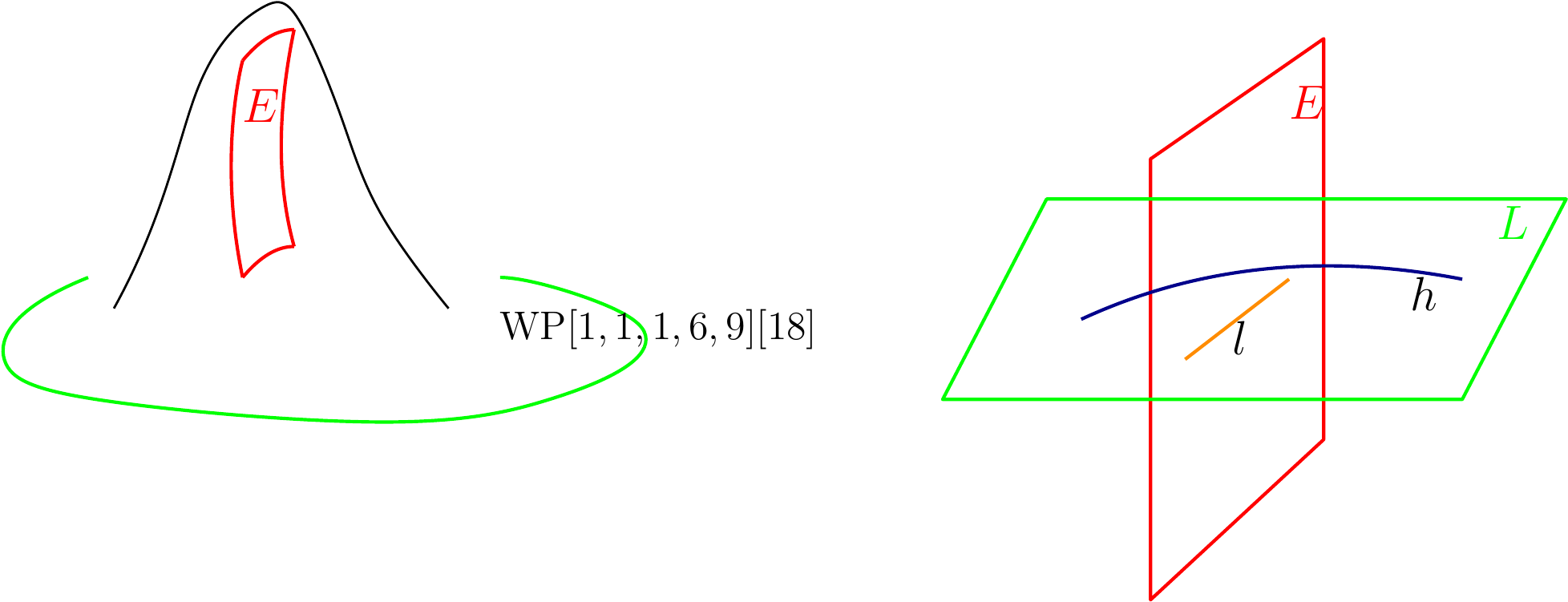}
 \caption{Elliptic f\/ibration in the A-model as a blow-up.} \label{figureblowup}
\end{figure}

Thinking of $\check{X}$ as the blow-up, the K\"ahler classes are linear combinations of $L$ (the strictly transform of the K\"ahler class on the singular variety) and the exceptional divisor class~$E$. The f\/ibration structure also tells that the class of the base~$\mathbb{P}^{2}$ is~$E$, the pull back of $\mathcal{O}_{\mathbb{P}^{2}}(1)$ gives the class~$L$. One has $E\cdot E\cdot L=-3$. This is the degree of the line bundle~$K_{\mathbb{P}^{2}}$ over $\mathbb{P}^{2}$. It conf\/irms the statement that~$E$ is the base~$\mathbb{P}^{2}$ of the elliptic curve family. The ef\/fective curve classes are
\begin{gather*}
h=L\cdot L,\qquad \ell=L\cdot E,
\end{gather*}
which represent the elliptic curve f\/iber and the hypersurface class in the base $E$, respectively. See the illustration in Fig.~\ref{figureblowup}. The dual nef cone is worked out to be the one generated by~$H$,~$L$.

These intersections are nicely encoded into the linear relations among the rays in the toric fan that def\/ines the ambient space
\begin{gather*}
Q_{1}=(1,1,1,-3,0,0;0),\qquad Q_{2}=(0,0,0,1,2,3;-6).
\end{gather*}
They represent the curve classes $ \ell$, $h$ respectively. The toric invariant divisors correspond to the columns of the above matrix of linear relations.
More precisely, one has
\begin{gather*}
L\sim \begin{pmatrix}1\\ 0 \end{pmatrix},\qquad H\sim\begin{pmatrix} 0\\ 1 \end{pmatrix},\qquad E\sim\begin{pmatrix}-3\\ 1 \end{pmatrix}.
\end{gather*}
The last vector $(0,-6)$ represents the f\/irst Chern class of $K_{\mathrm{WP[1,2,3]}}$, which is canonical sheaf of the f\/iber weighted projective space (in which the elliptic curve f\/iber sits). We denote the corresponding class by
\begin{gather*}
J\sim \begin{pmatrix} 0\\ -6 \end{pmatrix}.
\end{gather*}

Now a K\"ahler class is represented by a linear combination of these classes, with certain positivity conditions satisf\/ied,
\begin{gather*}
K=\log a_{1}L+\log a_{2}L+\log a_{3}L+\log a_{0}E+\log a_{4}(2H)+\log a_{5}(3H)+\log b_{0}J.
\end{gather*}
The parameters $a_{i}$, $b_{i}$ are mirror to the coordinates with the same names in the B-model, up to terms which do not af\/fect the qualitative analysis. In the following we shall use the same coordinates $a$, $b$ as in~\eqref{eqncomplexmoduli}. An element in the nef cone (one of the chambers in the second fan of the toric variety\footnote{See, e.g., \cite{Cox:1999} for a detailed review.}) must satisfy the condition
\begin{gather}
K=\log\left( {a_{1}a_{2}a_{3}\over a_{0}^{3}}\right)L +\log \left({a_{4}^{2}a_{5}^{3}a_{0}\over b_{0}^{6}}\right)H\nonumber\\
\hphantom{K}{} =\log \big(a^{-3}\big)L+\log \big(a b^{-6}\big) H \in \mathbb{R}_{>0}L\oplus \mathbb{R}_{>0}H.\label{eqnkahlerclass}
\end{gather}
Therefore the large volume limit corresponds to the point
\begin{gather}\label{eqnLCSL}
a^{-3}=0=ab^{-6},
\end{gather}
which is mirror to the large complex structure limit mentioned before. Hence one can see that in the def\/ining equation~$\xi$ in~\eqref{eqnmirrorfamily}
for the B-model, one needs to send both $a$, $b$ to $\infty$.

Now it is easy to see that the limit $b_{0}=0$ corresponds to the degeneration\footnote{However, in this limit $\log (ab^{-6})$ has a negative sign from the previous limit in~\eqref{eqnLCSL} and the description in~\eqref{eqnkahlerclass} for the K\"ahler class fails. This means that this limit does not sit inside the nef cone of $\check{X}$ and it has moved deeply into some other chamber in the secondary fan. This situation is typical in the so-called phase transition process, see~\cite{Witten:1993}.} of K\"ahler structure $K\rightarrow L$, as when considering for example instanton expansions cycles with inf\/inity volumes are suppressed. In particular, the class of the elliptic curve f\/iber $h=L\cdot L$ does not survive the limit. On the level of toric data or Mori cone of curve classes, the vector~$Q_{2}$ representing the class~$h$ is invisible in the limit and hence what is left is the one~$Q_{1}$ which is exactly the toric data that def\/ines the geometry~$K_{\mathbb{P}^{2}}$. On the level of geometry, that~$h$ has inf\/inite volume tells that the elliptic f\/ibration over $E$ gets decompactif\/ied to~$K_{\mathbb{P}^{2}}$. This is consistent with the B-model picture discussed below~\eqref{eqntoruscooordinatesdegeneration}.

Similarly, for the limit $a_{4}=a_{5}=0$, one has the degeneration $K\rightarrow L$. Now the curve class~$3\ell$, which is the one underlying the cubic in $E\cong \mathbb{P}^{2}$ survives this limit. This is also consistent with the B-model picture discussed below~\eqref{eqnGKZlimit}.

\section{Discussions and speculations}\label{secdiscussions}

For a general Laurent polynomial $F$, the relations among the monomials are conveniently described combinatorially by the Newton polytope. This in particular gives a shortcut in deriving the dif\/ferential equations.

As mentioned in Section \ref{secdifferencesofchainintegrals}, the dif\/ference between the universal family of cubics and the Hesse pencil is that the latter carries an extra level structure. This is what picks out the~$4$ monomials appearing in the Hesse pencil among the all the~$10$ cubic monomials. These~$4$ monomials is what determines the symmetries of the family and hence the GKZ symmetries. The linear space spanned by them encode the full data of the family.

Instead of using the Veronese embedding of $\mathbb{P}^{2}$ into $\mathbb{P}\check{H}^{0}(\mathbb{P}^2, \mathcal{O}(3))$, we use only the~$4$ monomials
$xyz$, $x^3$, $y^3$, $z^3$, $i=1,2,3$. Consider the partial Veronese map
\begin{gather*}
\Phi\colon \ \mathbb{P}^{2}\rightarrow \mathbb{P}^{3}, \qquad (x,y,z)\mapsto (X_{0}, X_{1},X_{2},X_{3})= \big(xyz,x^3,y^3,z^3\big).
\end{gather*}
The image of $\mathbb{P}^{2}$ is given by the vanishing of the ideal generated by
\begin{gather*}
\Phi=X_{1}X_{2}X_{3}-X_{0}^3,
\end{gather*}
which def\/ines a singular cubic surface $S_{0}$ of degree $3$ in~$\mathbb{P}^{3}$. Here we used the same notation $\Phi$ to denote the polynomial by abuse of notation. The singular points on $S_{0}$ are given by
\begin{gather*}
(X_{0},X_{1},X_{2},X_{3})=(0,1,0,0), \,(0,0,1,0),\, (0,0,0,1).
\end{gather*}
Locally a neighborhood of each singular point is of the form of an $A_{2}$-singularity described by~$\mathbb{C}^2/\mathbb{Z}_{3}$. This surface is nothing but the mirror of $\mathbb{P}^{2}$ given by $\mathbb{P}^{2}/G$ as described in~\eqref{eqnmirrorP2}. The elliptic curve is then mapped to the intersection of this cubic surface $S_{0}$ with the hyperplane
\begin{gather*}
H:= \sum a_{i} X_{i}=0.
\end{gather*}
One can alternatively make a moduli dependent embedding by using the monomials $a_{i}X_{i}$, \smash{$i=0,1,2,3$}, then the moduli dependence of the intersection is full encoded in the moduli dependence of the singular cubic surface~$S_{0}$.

The above procedure linearizes the monomials in the def\/ining equation $F$ for the Hesse pencil. The integrand in~\eqref{eqnosscilatingintegralaffine} then becomes (up to the $a_{0}$ factor)
\begin{gather*}
e^{-F}dxdydz={1\over 3^3} e^{-\sum\limits_{i=0}^{3} a_{i}X_{i}}\delta_{S_{0}}\prod {dX_{i}\over X_{i}^{2/3}}:=e^{-H}\delta_{S_{0}}\mu,
\end{gather*}
where $\delta$ means the Dirac delta distribution. As before, the integrand $e^{-sF}\mu_{0}\wedge ds$ on needs to be suitably interpreted in order to incorporate all the $\mathcal{Z}$-symmetries.

Now the original $\mathcal{Z}$-symmetries become diagonal symmetries of the quadratic form $H=\sum a_{i}X_{i}$. These symmetries also leave the variety $S_{0}$ invariant (the polynomial~$\Phi$ is not invariant). The original $\mathcal{D}$-symmetries are manifest through the Veronese map~$\Phi$. Schematically one has
\begin{gather}
\mathcal{Z}(H)=0,\qquad\! \mathcal{Z}(S_{0})=0, \qquad\! \mathcal{Z}\big(e^{-H}\delta_{S_{0}}\mu\big)=0,\qquad\!
\mathcal{D}\big(e^{-H}\delta_{S_{0}}\mu\big)=e^{-H}\delta_{S_{0}} \Phi \mu=0.\!\!\!\!\label{eqnDmodule}
\end{gather}

Note that the characteristic variety (singular support) of $\mathcal{D}$ is now def\/ined by
\begin{gather*}
\big\{p_{1}p_{2}p_{3}-p_{0}^{3}=0\big\}\subseteq T^{*}\mathbb{C}^{4},
\end{gather*}
where $p_{i}$, $i=0,1,2,3$ are the f\/iber coordinates of the cotangent bundle of $\mathbb{C}^{4}$ parametrized by $a_{0}$, $a_{1}$, $a_{2}$, $a_{3}$.
This takes the same form as the polynomial~$\Phi$. In fact, by restricting the f\/ield~$\mathbb{C}$ to the f\/ield~$\mathbb{R}$ which does not af\/fect the discussion on GKZ symmetries, one can regard
\begin{gather*}
\mathcal{F}(\bullet)=\int e^{-H} (\bullet)\mu
\end{gather*}
as a formal Fourier transform. Then one has
\begin{gather*}
\mathcal{D}=\mathcal{F}(\Phi) .
\end{gather*}
Note that $\mathcal{D}$ has constant coef\/f\/icients, this then relates $\Phi$ to the characteristic variety.

We can also formally write the original oscillating integral over the invariant chain $D_{3}=(0,\infty)\times (0,\infty)\times (0,\infty)$ in \eqref{eqn3dcontour} (which produces all solutions under the monodromy action) as
\begin{gather*}
I=\int_{D_{3}}e^{-F}\mu_{0}=\int_{D_{3}\times (0,\infty)}e^{-H}\delta_{S_{0}}\mu=\mathcal{F}(\delta_{S_{0}}).
\end{gather*}

The conf\/iguration $(\Phi,H)$ fully encodes the information of the Hesse pencil. Since by a moduli dependent embedding one can make $H$ independent of $a_{1}$, $a_{2}$, $a_{3}$, $a_{0}$, it is therefore very natural to expect that the GKZ system can be approached via studying the mixed Hodge structure of $S_{0}$ together with the moduli-independent hyperplane. The presentation of the GKZ symmetries in the form displayed in~\eqref{eqnDmodule} also begs for an explanation in terms of D-module in addressing this problem. A good understanding of this matter is potentially useful in the studies of open mirror symmetry where similar situation occurs, see, e.g.,~\cite{Lian:2012, Mayr:2001}.

\subsection*{Acknowledgements}

The author dedicates this article to Professor Noriko Yui on the occasion of her birthday. The author is grateful for her constant encouragement and support, and in particular for many inspiring discussions on geometry and number theory. The author would like to thank Murad Alim, An Huang, Bong Lian and Shing-Tung Yau for discussions on open string mirror symmetry which to a large extent inspired this project. He thanks further Kevin Costello, Shinobu Hosono, Si Li and Zhengyu Zong for their interest and helpful conversations on Landau--Ginzburg models and chain integrals, and Don Zagier for some useful discussions on modular forms back in year 2013.
He also thanks the anonymous referees whose suggestions have helped improving the article.

This research was supported in part by Perimeter Institute for Theoretical Physics. Research at Perimeter Institute is supported by the Government of Canada through Innovation, Science and Economic Development Canada and by the Province of Ontario through the Ministry of Research, Innovation and Science.

\pdfbookmark[1]{References}{ref}
\LastPageEnding

\end{document}